\documentclass[a4paper,10pt]{amsart}
\usepackage{amssymb,amsmath}
\usepackage{amsfonts}
\usepackage[T1]{fontenc}

\newtheorem{theorem}{Theorem}[section]
\newtheorem{lemma}[theorem]{Lemma}
\newtheorem{proposition}[theorem]{Proposition}
\newtheorem{corollary}[theorem]{Corollary}
\newtheorem{claim}{Claim}
\newtheorem{definition}[theorem]{Definition}
\newtheorem*{notation}{Notation}

\theoremstyle{remark}
\newtheorem{remark}[theorem]{Remark}

\newcommand{\bl}{\begin{lemma}}
\newcommand{\el}{\end{lemma}}
\newcommand{\bpr}{\begin{proposition}}
\newcommand{\epr}{\end{proposition}}
\newcommand{\bp}{\begin{proof}}
\newcommand{\ep}{\end{proof}}
\newcommand{\bd}{\begin{definition}}
\newcommand{\ed}{\end{definition}}
\newcommand{\bt}{\begin{theorem}}
\newcommand{\et}{\end{theorem}}
\newcommand{\bc}{\begin{corollary}}
\newcommand{\ec}{\end{corollary}}
\newcommand{\br}{\begin{remark}}
\newcommand{\er}{\end{remark}}
\newcommand{\bcl}{\begin{claim}}
\newcommand{\ecl}{\end{claim}}

\newcommand{\N}{{\mathbb{N}}}
\newcommand{\abs}[1]{\lvert#1\rvert}
\newcommand{\nrm}[1]{\|#1\|}
\newcommand{\al}{\alpha}
\newcommand{\e}{\varepsilon}
\newcommand{\de}{\delta}
\newcommand{\bnum}{\begin{enumerate}}
\newcommand{\enum}{\end{enumerate}}
\newcommand{\mc}{\mathcal}

\numberwithin{subsection}{section}
\numberwithin{equation}{section}

\DeclareMathOperator{\supp}{supp}
\DeclareMathOperator{\maxsupp}{maxsupp}
\DeclareMathOperator{\minsupp}{minsupp}
\DeclareMathOperator{\ran}{range}
\DeclareMathOperator{\suc}{succ}

\begin{document}
\title{Quasiminimality in mixed Tsirelson spaces}
\author{Antonis Manoussakis and Anna Maria Pelczar}
\address{Department of Mathematics, University of Aegean, Karlovasi,
Samos, GR 83200, Greece} \email{amanouss@aegean.gr}
\address{Institute of Mathematics, Jagiellonian University, {\L}ojasiewicza 6, 30-348 Krak\'ow, Poland}
\email{anna.pelczar@im.uj.edu.pl}
\begin{abstract}
We prove quasiminimality of the regular mixed Tsirelson spaces
$T[(\mathcal{S}_{n},\theta_{n})_{n}]$ with the sequence
$(\frac{\theta_n}{\theta^n})_n$ decreasing, where $\theta=\lim_n
\theta_n^{1/n}$, and quasiminimality of all mixed Tsirelson spaces
$T[(\mathcal{A}_{n},\theta_{n})_{n}]$. We prove that under certain
assumptions on the sequence $(\theta_n)_{n}$ the dual spaces are
quasiminimal.
\end{abstract}
\maketitle

\section*{Introduction}

Mixed Tsirelson spaces form an important class of spaces in the theory of Banach spaces
extensively studied with respect to their distortability, local structure as well as
minimality properties.

Recall that an infinite dimensional Banach space $X$ is
\textit{minimal}, if any closed infinite dimensional subspace of
$X$ contains a further subspace isomorphic to $X$. A infinite
dimensional Banach space $X$ is \textit{quasiminimal}, if any two
infinite dimensional subspaces $Y,Z$ of $X$ contains further two
infinite dimensional subspaces $Y'\subset Y, Z'\subset Z$ which
are isomorphic. The use of Ramsey theory in the famous Gowers'
dichotomy inspired studies on minimality and related properties,
e.g. \cite{fr,jko,dfko}, in particular in mixed Tsirelson spaces,
cf. also \cite{klmt,lt,m}.

We recall here briefly results concerning minimality properties of mixed Tsirelson
spaces. The famous Schlumprecht space $S=T[(\mc{A}_n,1/\log_2(n+1))_n]$ is complementably
minimal by \cite{as}, i.e. any closed infinite dimensional subspace of $S$ contains a
complemented copy of the whole space. This result was proved also for a class of
superreflexive spaces extending the construction of $S$ in \cite{CKKM} and for certain
types of mixed Tsirelson spaces $T[(\mc{A}_{k_n},\theta_n)_n]$ in \cite{m}. On other hand Tzafriri space $T[(\mc{A}_n,\frac{1}{2\sqrt{n}})_n]$ is not minimal \cite{jko}.

In the original Tsirelson space $T[\mc{S}_1,1/2]$ every normalized
block sequence is equivalent to a subsequence of the basis. In
context of mixed Tsirelson spaces $T[(\mc{S}_{k_n},\theta_n)_n]$ a
more general property - subsequential minimality - was studied. A
space with a basis is called is \textit{subsequentially minimal},
if any block subspace contains a normalized block sequence
equivalent to a subsequence of the basis. In \cite{lt} it was
proved that if a regular sequence $(\theta_n)$ satisfies the
condition $(\star)$:
$$
\lim_m\limsup_n\frac{\theta_{m+n}}{\theta_n}>0
$$
then the mixed Tsirelson space $X=T[(\mc{S}_n,\theta_n)_n]$ is subsequentially minimal if
and only if any block subspace of $X$ admits an $\ell_1-\mc{S}_\omega-$spreading model,
if and only if any block subspace of $X$ has Bourgain $\ell_1-$index greater than
$\omega^\omega$. In particular $X$ is complementably subsequentially minimal if
$\sup\theta_n^{1/n}=1$ \cite{m}, i.e. subspaces shown to be isomorphic to subspaces
generated by subsequences of the basis are complemented in the space. In \cite{klmt}
analogues of these results were studied in the partly modified setting, in particular it
was shown that a partly modified mixed Tsirelson space is subsequentially minimal
provided $\sup\theta_n^{1/n}=1$. Also in \cite{klmt} a large class of mixed Tsirelson
spaces failing the subsequential minimality in a strong sense was described.

In the paper we study $p-$spaces $T[(\mc{A}_n,\theta_n)_n]$, as
defined in \cite{m}, as well as mixed Tsirelson spaces
$T[(\mc{S}_n,\theta_n)_n]$ showing quasiminimality of all spaces
of the first type and regular spaces of the second type provided
the sequence $(\frac{\theta_n}{\theta^n})_n$ is decreasing, where
$\theta=\lim_n \theta_n^{1/n}$, and we consider related
properties. Notice that for the spaces of the second type, in case
$\theta<1$, the monotonicity of the sequence
$(\frac{\theta_n}{\theta^n})_n$ excludes the condition $(\star)$
of \cite{lt} mentioned above.

In the study of spaces of both types special averages provide the
major tool. In case of $p-$space we show that $\ell_p$ is finitely
block represented in all block subspaces of any $p-$spaces, i.e.
$p$ is in Krivine set of all block subspaces of a $p-$space.
Equivalent block sequences in these spaces are formed by long
averages of $\ell_p-$averages.

In case of mixed Tsirelson spaces of type
$T[(\mc{S}_n,\theta_n)_n]$ we use special averages, whose
properties were studied in \cite{ao}, vectors of the same type as
special convex combinations used in \cite{ad,adm}. In this case,
thanks to the higher complexity of the families $(\mc{S}_n)_n$,
equivalent sequences are formed by long special averages (not
averages of special averages, as it is in $p-$spaces case). The
reasoning proving the quasiminimality of spaces
$T[(\mc{S}_n,\theta_n)_n]$ implies also the result of \cite{m} on
subsequential minimality of the spaces in case
$\sup\theta_n^{1/n}=1$. In both cases the presented argument gives
saturation of considered mixed Tsirelson spaces
$T[(\mc{S}_n,\theta_n)_n]$ or $T[(\mc{A}_n,\theta_n)_n]$ with
subsequentially minimal subspaces.

Notice that in both classes of mixed Tsirelson spaces special averages were used in the
proofs of arbitrary distortability of spaces under additional assumptions. The
observations presented in the paper indicates also similarities in the way spaces
$\ell_p$, $1\leq p<\infty$, and Tsirelson type spaces $T[\mc{S}_1,\theta]$, $0<\theta\leq
1$, are represented in spaces $T[(\mc{A}_n,\theta_n)_n]$ and $T[(\mc{S}_n,\theta_n)_n]$.

Now we describe the contents of the paper. In Section 1 we recall basic notions. In
Section 2. we consider $p-$spaces, state technical lemmas needed in the sequel and
describe representation of $\ell_p$ in $p-$spaces $T[(\mc{A}_n,\theta_n)_n]$ depending on
the behavior of $(\theta_n)_n$. The essential tool for quasiminimality in $p-$spaces case
is provided by Theorem \ref{kriv} stating that the Krivine set of any subspace of a
$p-$space, $1\leq p<\infty$ contains $p$. Recall here that by \cite{m} Krivine set of a
$p-$space $T[(\mc{A}_n,\theta_n)_n]$ contains only $p$ provided the sequence
$(\theta_nn^{1/q})_n$ is decreasing, for $1/p+1/q=1$, while there are $p-$spaces, $p>1$,
with 1 in the Krivine set. We present also analogues to observations concerning the
original Tsirelson space and mixed Tsirelson spaces given in \cite{ot,otw}, in particular
extending the result of \cite{jko} on saturation of Tzafriri space by $\ell_2-$asymptotic
subspaces. Section 3 contains the proof of quasiminimality of arbitrary $p-$spaces
(Theorem \ref{p-quasi}) and study of behavior of $p-$spaces, $1<p<\infty$, and their
duals in case $\inf_n\theta_nn^{1/q}>0$ (Theorem \ref{tz-quasi}), including in particular
Tzafriri space. Section 4 is devoted to mixed Tsirelson spaces
$T[(\mc{S}_n,\theta_n)_n]$. After technical lemmas concerning Schreier families, we show
the existence of special averages (Corollary \ref{t-av}), applying the reasoning of
\cite{ao} in case the sequence $(\frac{\theta_n}{\theta^n})_n$ is decreasing. Under this
assumption, using special averages we prove the quasiminimality of a regular mixed
Tsirelson spaces $T[(\mc{S}_n,\theta_n)_n]$ (Theorem \ref{t-quasi}). Section 5 presents a
general argument enabling us to transfer the minimality and quasiminimality properties of
spaces to their duals provided certain complementability property holds. In particular we
prove quasiminimality of duals to mixed Tsirelson spaces $T[(\mc{S}_n,\theta_n)_n]$ in
case $\sup\theta_n^{1/n}=1$.

\section{Preliminaries}
We recall now the basic definitions and standard notation. Let $X$
be a Banach space with a basis $(e_i)$. The \textit{support} of a
vector $x=\sum_i x_i e_i$ is the set $\supp x =\{ i\in \N :
x_i\neq 0\}$, the \textit{range} of a vector $x\in X$ - the
smallest interval in $\N$ containing support of $x$. Given any
$x=\sum_i x_ie_i$ and finite $E\subset\N$ put $Ex=x_E=\sum_{i\in
E}x_ia_i$. We write $x<y$ for vectors $x,y\in X$, if $\max\supp
x<\min \supp y$. A \textit{block sequence} is any sequence
$(x_i)\subset X$ satisfying $x_{1}<x_{2}<\dots$, a \textit{block
subspace} of $X$ - any closed subspace spanned by an infinite
block sequence. A subspace spanned by a basic sequence $(x_n)$ we
denote by $[x_n]$.

A Banach space $X$ is \textit{saturated} with subspaces of type
$(*)$, if any infinite dimensional subspace of $X$ contains a
further infinite dimensional subspace of type $(*)$.

A basic sequence $(x_n)$ $C-$\textit{dominates} a basic sequence
$(y_n)$, $C\geq 1$, if for any scalars $(a_n)$ we have
$$
\nrm{\sum_na_ny_n}\leq C\nrm{\sum_na_nx_n}
$$
Two basic sequences $(x_n)$ and $(y_n)$ are
$C$-\textit{equivalent}, $C\geq 1$, if $(x_n)$ $C-$dominates
$(y_n)$ and $(y_n)$ $C-$dominates $(x_n)$.

Let $\mc{M}$ be a family of finite subsets of $\N$. We say that $\mc{M}$ is
\textit{regular}, if it is \textit{hereditary}, i.e. for any $G\subset F$, $F\in \mc{M}$
also $G\in \mc{M}$, \textit{spreading}, i.e. for any integers $n_1<\dots<n_k$ and
$m_1<\dots<m_k$ with $n_i\leq m_i$, $i=1,\dots, k$, if $(n_1,\dots,n_k)\in \mc{M}$ then
also $(m_1,\dots,m_k)\in \mc{M}$, and \textit{compact} in the product topology of $2^\N$.
Given a spreading family $\mc{M}$ we say that a sequence $E_1<\cdots<E_n$ of finite
subsets of $\N$ is $\mc{M}$-\textit{admissible}, if $(\min E_i)_{i=1}^n\in\mc{M}$. For
any two families $\mc{M},\mc{N}$ of finite subsets of $\N$ we define
$$
\mc{M}[\mc{N}]=\{F_1\cup\dots\cup F_k:\ F_1,\dots,F_k\in\mc{N}, F_1<\dots<F_k -
\mc{M}-admissible\}
$$

We will work on two types of regular families:
$(\mc{A}_n)_{n\in\N}$ and $(\mc{S}_\al)_{\al<\omega_1}$. Let
$$
\mc{A}_n=\{F\subset\N:\# F\leq n\}, \ \ n\in\N
$$
\textit{Schreier families} $(\mc{S}_\al)_{\al<\omega_1}$,
introduced in \cite{aa}, are defined by induction:
$$
\mc{S}_0 =\{\{ n\}:\ n\in\N\}\cup\{\emptyset\}, \ \ \mc{S}_{\al+1} =\mc{S}_1[\mc{S}_\al], \ \
\al<\omega_1
$$
If $\al$ is a limit ordinal, choose $\al_n\nearrow \al$ and set
$$
\mc{S}_\al=\{F:\ F\in \mc{S}_{\al_n}\ \mathrm{and}\ n\leq F\ \mathrm{for\ some}\ n\in\N\}
$$
It is well known that all families $(\mc{A}_n)_{n\in\N}$, $(\mc{S}_\al)_{\al<\omega_1}$ are regular.

Let $X$ be a Banach
space with a basis. We say that a sequence $x_1<\dots <x_n$ is
$\mc{M}$-\textit{admissible}, if $(\supp x_i)_{i=1}^n$ is $\mc{M}$-admissible.

\bd[Tsirelson-type space] Fix a regular family $\mc{M}$ and $0<\theta<1$. Let $K\subset
c_{00}$ be the smallest set satisfying the following:

\bnum \item $(\pm e_n^*)_n\subset K$,

\item for any functionals $\phi_1<\dots<\phi_k$ in $K$, if
$(\phi_i)_{i=1}^k$ is $\mc{M}$-admissible, then
$\theta(\phi_1+\dots+\phi_k)\in K$.

\enum We define a norm on $c_{00}$ by the formula
$\nrm{x}=\sup\{f(x):f\in K\}$, $x\in c_{00}$. Then the
\textit{Tsirelson-type space} $T[\mc{M},\theta]$ is the completion
of $(c_{00}, \nrm{\cdot})$. \ed It is standard to verify that the
norm $\nrm{\cdot}$ is the unique norm on $c_{00}$ satisfying the
equation
$$
\nrm{x}=\max\left\{\nrm{x}_\infty,\sup\left\{\theta\sum_{i=1}^k\nrm{E_ix}: \ (E_i)_{i=1}^k -
\mc{M}- \mathrm{admissible}\right\}\right\}
$$
Recall that $T[\mc{A}_n,\theta]=c_0$ if $\theta\leq 1/n$ and $T[\mc{A}_n,\theta]=\ell_p$,
if $\theta=1/\sqrt[q]{n}$ for $q$ satisfying $1/p+1/q=1$. The space $T[\mc{S}_1,1/2]$ is
the Tsirelson space, the first discovered non-classical space - not containing $\ell_p$,
$1\leq p<\infty$ or $c_0$.

Now we recall the definition of the mixed Tsirelson spaces, using
not only one family $\mc{M}$ and parameter $\theta$, but sequences
of these objects. Again the norm can be defined by its norming set
or as a solution of an equation.

\bd[Mixed Tsirelson space] Fix a sequence of regular families $(\mc{M}_n)$ and sequence
$(\theta_n)\subset (0,1)$ with $\lim_{n\to\infty}\theta_n=0$. Let $K\subset c_{00}$ be
the smallest set satisfying the following:

\bnum \item $(\pm e_n^*)_n\subset K$,

\item for any $\phi_1<\dots<\phi_k$ in $K$, if
$(\phi_i)_{i=1}^k$ is $\mc{M}_n$-admissible for some $n\in\N$,
then $\theta_n(\phi_1+\dots+\phi_k)\in K$.

\enum We define a norm on $c_{00}$ by $\nrm{x}=\sup\{f(x):f\in
K\}$, $x\in c_{00}$. The \textit{mixed Tsirelson space}
$T[(\mc{M}_n,\theta_n)_n]$ is the completion of $(c_{00},
\nrm{\cdot})$. \ed

Again it is standard to verify that the norm $\nrm{\cdot}$  is the unique
norm on $c_{00}$ satisfying the equation
$$
\nrm{x}=\max\left\{\nrm{x}_\infty,\sup\left\{\theta_n\sum_{i=1}^k\nrm{E_ix}: \ (E_i)_{i=1}^k -
\mc{M}_n- \mathrm{admissible}, \ n\in\N\right\}\right\}
$$
It follows immediately that the unit vector basis $(e_n)$ is 1-unconditional in the space
$T[(\mc{M}_n,\theta_n)_n]$. It was proved in \cite{ad} that any
$T[(\mc{S}_{k_n},\theta_n)_n]$ is reflexive, also any $T[(\mc{A}_{k_n},\theta_n)_n]$ is
reflexive, provided $\theta_n>\frac{1}{k_n}$ for at least one $n\in\N$. In this setting
Schlumprecht space $S$ is the space $T[(\mc{A}_n,\frac{1}{\log_2(n+1)})_n]$, Tzafriri
space is $T[(\mc{A}_n,\frac{c}{\sqrt{n}})_n]$ for $0<c<1$.

A \textit{tree} is a partially ordered set $(T,\prec)$ such that for any $t\in T$ the set
$\{s\in T: \ s\prec t\}$ is finite and linearly ordered. Given a tree $T$ for any $t\in
T$ by $\suc t$ we denote the set of immediate successors of $t$ in $T$. The
\textit{height} of a tree $T$ is the supremum of the length of branches in the tree $T$
(i.e. well ordered subsets of the tree $T$). The \textit{level} of an element $t$ of a
tree $T$ with a unique root (i.e. the minimal element) is the length of the branch
linking the root with $t$.

The following notion provides a useful tool for estimating norms in mixed Tsirelson
spaces:

\bd[The tree-analysis of a norming functional] Let $\phi\in K$. By a
\textit{tree-analysis} of $\phi$ we mean a finite family $(\phi_t)_{t\in T}$ indexed by a
tree $T$ with a unique root $0\in T$ such that the following hold

\bnum \item $\phi_{0}=\phi$ and $\phi_t\in K$ for all $t\in T$, \item $t\in T$ is maximal
if and only if $\phi_t\in (\pm e_n^*)$, \item for every $t\in T$ not maximal there is
some $n\in\N$ such that $(\phi_s)_{s\in\suc t}$ is an $\mc{M}_n$-admissible block
sequence and $\phi_t=\theta_n(\sum_{s\in\suc t}\phi_s)$. We call $\theta_n$ the
\textit{weight} of $\phi_t$. \enum \ed Notice that every functional from a norming set
$K$ admits a tree-analysis, not necessarily unique.
\begin{definition} Let $(y_n)\subset X$ be a block sequence.

a)  Let $\phi\in K$ be a  norming functional  with a
tree-analysis $(\phi_t)_{t\in T}$. We say that $\phi_t$ \textit{covers} $y_n$ if $t\in T$
is maximal in $T$ with the property $\supp \phi_t\supset \supp y_n\cap \supp \phi$.

b) a finite sequence  $(E_{k})_{k=1}^{m}$  of subsets of $\N$ is
said to be \textit{comparable} with the sequence $(y_n)$ if for
every $n$ and every $k\leq m$ we have
$$
E_{k}\subset \ran y_{n}\,\,\, or\,\,\, \supp
y_{n}\cap\bigcup_{i=1}^mE_i\subset E_{k} \,\,\,or\,\,\,\ran
E_{k}\cap \ran y_{n}=\emptyset.
$$
c) A  norming functional $\phi\in K$ with a tree-analysis
$(\phi_t)_{t\in T}$ is said to be \textit{comparable} with
$(y_n)_{n}$ if the sequence $(\supp\phi_t)_{t\in T}$ is comparable
with $(y_{n})_{n}$.
\end{definition}
\section{Representability of $\ell_p$ in $p-$spaces}
We describe now the class of $p-$spaces following the definition of \cite{m}. Recall that
any space $T[(\mc{A}_{k_n},\theta_n)_{n=1}^N]$, defined analogously to Tsirelson type
spaces and mixed Tsirelson spaces, with the use of finitely many families
$(\mc{A}_{k_n})_{n=1}^N$ and finitely many scalars $(\theta_n)_{n=1}^N\subset (0,1)$, is
isomorphic to some $\ell_p$, $1\leq p <\infty$ or $c_0$ \cite{bd}.

\bd \cite{m} A mixed Tsirelson space $T[(\mc{A}_{k_n},\theta_n)_{n\in\N}]$ is called a
\textit{$p$-space}, for $p\in [1,\infty)$, if there is a sequence $(p_N)_N\subset
(1,\infty)$ such that

\bnum \item $p_N\to p$ as $N\to\infty$, and $p_N\geq p_{N+1}>p$ for any $N\in\N$, \item
$T[(\mc{A}_{k_n},\theta_n)_{n=1}^N]$ is isomorphic to $\ell_{p_N}$ for any $N\in\N$.
\enum

A $p-$space $T[(\mc{A}_n,\theta_n)_{n\in\N}]$ is called \textit{regular}, if
$\theta_n\searrow 0$ and $\theta_{nm}\geq\theta_n\theta_m$ for any $n,m\in\N$.

\ed

Notice that any mixed Tsirelson space $T[(\mc{A}_{k_n},\theta_n)_{n\in\N}]$ is isometric
to a space $T[(\mc{A}_n,\widehat{\theta}_n)_{n\in\N}]$, where
$$
\widehat{\theta}_n=\sup\left\{\prod_{i=1}^l\theta_{n_i}:\ \prod_{i=1}^l n_i\geq n
\right\}, \ \ \ n\in\N
$$
It follows directly that if $T[(\mc{A}_{k_n},\theta_n)_{n\in\N}]$ is a $p-$space, then
also $T[(\mc{A}_n,\widehat{\theta}_n)_{n\in\N}]$ is a $p-$space. Notice also that the
sequence $(\widehat{\theta}_n)$ satisfies $\widehat{\theta}_n\searrow 0$,
$\widehat{\theta}_{nm}\geq \widehat{\theta}_n\widehat{\theta}_m$ for any $n,m\in\N$.

Therefore from now on we assume that $p-$spaces we consider are regular.

\begin{notation}

Let $T[(\mc{A}_n,\theta_n)_{n\in\N}]$ be a regular $p-$space. If we set
$\theta_n=1/n^{1/q_n}$ with $q_n\in (1,\infty)$, $n\in\N$, then $q=\lim_n q_n=\sup_n
q_n\in (0,\infty]$, where $1/p+1/q=1$, with usual convention $1/\infty=0$.

In the situation as above let $c_n=\theta_nn^{1/q}\in (0,1)$, $n\in\N$, if $p>1$. To
unify the notation put $c_n=\theta_n$, $n\in\N$, in case $p=1$.
\end{notation}
\br In \cite{m} in the definition in the $p-$spaces it was posed
also that the sequence $(p_n)_{n}$ was strictly decreasing and
hence $\theta_{n}$ also. It was also posed that $c_n\searrow 0$.
In our setting, the sequences $(p_n)_{n}, (c_{n})_{n}$ do not
necessarily satisfy these two conditions.
\er

\br\label{rem-pnorm} Notice that the norming set of a $p-$space is
contained in the unit ball of $\ell_q$, hence
$\nrm{\cdot}\leq\nrm{\cdot}_p$.
\er

First we present two technical observations.

\bl \label{p-xi} Let $X$ be a $p-$space, $1\leq p<\infty$. Take a functional $\phi$ with a tree-analysis $(\phi_t)_{t\in T}$ and
any finite block sequence $(v_n)$. Set $v=\sum_n v_n$. Then there is a functional $\phi'$
with a tree-analysis $(\phi'_t)_{t\in T'}$ comparable with $(v_n)$ satisfying
$6\phi'(v)\geq \phi (v)$. \el

\bp

We can assume that $\supp\phi\subset\bigcup_n\supp v_n$. For any $n$ pick $t_n\in T$  so
that $\phi_{t_n}$ covers $v_n$. Let $E_1<\dots<E_k$ be supports of all successors
$\phi_1,\dots,\phi_k$ of $\phi_{t_n}$ which intersect the support of $v_n$.

CASE 1. There is $i=1,\dots, k$ so that $E_i\subset \supp v_n$. If $E_1\cap \supp
v_{n-1}\neq\emptyset$ and $\phi_1(v_n)\geq \phi_i(v_n)$ then change the tree: split
$\phi_1$ into two parts supported on $E_1\cap \supp v_n$ and $E_1\setminus \supp v_n$ and
erase $\phi_i$, put the part $\phi_1|_{\supp v_n}$ into tree instead of $\phi_i$, leave
$\phi_1|_{E_1\setminus \supp v_n}$ in the place of $\phi_1$. If $E_1\cap \supp
v_{n-1}\neq\emptyset$ and $\phi_1(v_n)\leq \phi_i(v_n)$ then erase from $\phi$ the part
supported on $E_1\cap \supp v_n$. Analogously proceed, if $E_k\cap \supp
v_{n+1}\neq\emptyset$. The action of the modified functional on $v_n$ is not less than
$\phi(v_n)/3$.

CASE 2. $k=2$, $E_1\cap \supp v_{n-1}\neq\emptyset$, $E_2\cap
\supp v_{n+1}\neq\emptyset$. Let for example $\phi_1(v_n)\geq
\phi_2(v_n)$, then erase part of $\phi$ supported on $E_2\cap
\supp v_n$, and for the modified functional find new $t_n$ - then
we have only Case 1. \ep

\bl\label{p-est} Let $X$ be a $p-$space, $1\leq p<\infty$. Take any norming functional
$\phi$ with a tree-analysis $(\phi_t)_{t\in T}$ and some $J_1,\dots,J_N\subset T$ so that

\bnum

\item $J_n\subset \suc t_n$ for some $t_n\in T$ for any $n=1,\dots,N$,

\item $\supp \phi_i\cap \supp \phi_j=\emptyset$ for any $i\neq j$,
$i,j\in \bigcup_n J_n$

\item $\supp\phi=\bigcup\{\supp \phi_i:\ i\in \bigcup_n J_n\}$

\enum

Let $\phi=\sum_{n=1}^N\gamma_n\sum_{i\in J_n}\phi_i$. Then for any scalars
$a_1,\dots,a_N\geq 0$ we have
$$
\sum_{n=1}^N a_n\gamma_n(\# J_n)^{1/q}\leq (\sum_{n=1}^Na_n^p)^{1/p}
$$
\el

\bp The case $p=1$ is obvious. Assume $p>1$. We proceed by induction on the height of the
tree $T$. For the height equal 0 (i.e. $\phi$ is a unit vector) the result is clear.
Assume we have the result for functionals with the tree-analysis of height smaller than
$l$ for some $l\geq 1$ and pick a functional with the tree-analysis $(\phi_t)_{t\in T}$
of height $l$. Let $J_1,\dots,J_N$, $\gamma_1,\dots,\gamma_N$, $a_1,\dots,a_N$ be as in
the Lemma.

Let $\phi=c_m/m^{1/q}\sum_{t\in \suc 0}\phi_t$. Let $I=\{n=1,\dots,N:\ J_n\subset \suc
0\}$. Let $J=\suc 0\setminus\bigcup_n J_n$. For any $t\in J$ let $I_t=\{n\not\in I:\ i\succ t\
\forall i\in J_n\}$. By the inductive assumption for $\phi_t$, $t\in J$, and H\"older
inequality we have
\begin{align*}
\sum_{n=1}^N\gamma_na_n(\# J_n)^{1/q} &\leq \frac{1}{m^{1/q}}\sum_{n\in I}a_n(\#
J_n)^{1/q}+\frac{1}{m^{1/q}}\sum_{t\in J}\sum_{n\in I_t}a_n\gamma_n m^{1/q}(\# J_n)^{1/q}
\\ &\leq \frac{1}{m^{1/q}}(\sum_{n\in I}a_n^p)^{1/p}(\sum_{n\in I}\#
J_n)^{1/q}+\frac{1}{m^{1/q}}\sum_{t\in J}(\sum_{n\in I_t}a_n^p)^{1/p}
\\ &\leq\frac{1}{m^{1/q}}(\sum_{n\in I}a_n^p)^{1/p}(\sum_{n\in I}\#
J_n)^{1/q}+\frac{1}{m^{1/q}}(\# J)^{1/q}(\sum_{t\in J}\sum_{n\in I_t}a_n^p)^{1/p}
\end{align*}

Notice that $I\cup\bigcup_{t\in J}I_t=\{1,\dots,N\}$ and $\#
J+\sum_{n\in I}\# J_n= \# \suc 0\leq m$, hence again by H\"older
inequality we obtain the desired upper bound by
$(\sum_{n=1}^Na_n^p)^{1/p}$. \ep

We recall here the notion of  $\ell_r-$averages - the basic tool we will use in studying
the properties of $p-$spaces:

\bd A vector $x\in X$ is called a $C-\ell_r-$\textit{average} of length $m$, for $r\in
[1,\infty]$, $m\in\N$ and $C\geq 1$ if $x=\sum_{n=1}^mx_i/\nrm{\sum_{n=1}^mx_i}$ for some
normalized block sequence $(x_n)_{n=1}^m$ which is $C$-equivalent to the unit vector
basis of $\ell_r^m$.

\ed

The following lemma extends the standard result for $\ell_1-$averages.

\bl\label{p-av}

Fix $1\leq r <\infty$ and $M\in\N$. Let $N\geq (2M)^r$. Then for any $C-\ell_r-$average
$x\in X$ of length $N$ and any $j\leq M$ we have
$$
\frac{j^{1/s}}{2C^2}\leq\sup\left\{\sum_i\nrm{E_ix}:\ (E_i) \ - \
\mc{A}_j\ - \mathrm{admissible}\right\}\leq 2C^2j^{1/s}
$$
where $1/s+1/r=1$.

\el

\bp

Let $x=(x_1+\dots+x_N)/\nrm{x_1+\dots+x_N}$, for some normalized sequence $x_1<\dots<x_N$
$C$-equivalent to the unit vector basis of $\ell_p^N$. Fix $j\leq M$. For the lower
estimate take subsets $I_1<\dots<I_j$ of $\{1,\dots,N\}$ such that $\# I_i=[N/j]$, i.e.
the integer part of $N/j$, let $E_i=\bigcup_{n\in I_i}\supp x_n$ and compute
\begin{align*}
\sum_{i=1}^j\nrm{E_ix} \geq\sum_i\nrm{\sum_{n\in I_i}
x_n}/CN^{1/r}\geq
\frac{j}{C^2N^{1/r}}\left[\frac{N}{j}\right]^{1/r}\geq\frac{j^{1/s}}{2^{1/r}C^2}\geq
\frac{j^{1/s}}{2C^2}.
\end{align*}
Now take any $E_1<\dots<E_j$ for $j\leq M$. For any $i=1,\dots,j$ pick $I_i=\{n:\ \supp x_n\cap
E_i\neq\emptyset\}$ and $I_i'=\{n\in I_{i}: \supp x_n\cap E_l=\emptyset\,\,\forall i\ne l\}$ and compute
$$
\nrm{E_i(x_1+\dots+x_N)}\leq \nrm{x_{\min I_i}} +\nrm{x_{\max I_i}}+
\nrm{E_{i}(\sum_{n\in I_{i}'}x_n)} \leq 2+C(\#I_{i}')^{1/r}
$$
Hence summing over $i$ we get
$$
\sum_i\nrm{E_ix}\leq 2M\frac{C}{N^{1/r}}+C^2\sum_{i}\frac{(\#I_{i}')^{1/r}}{N^{1/r}}\leq
2C^2j^{1/s}
$$
\ep

Now we study the way $\ell_p$ is represented in a $p-$space. Recall that by Krivine
theorem for any Banach space $X$ with a basis there is some $1\leq p\leq \infty$ such
that $\ell_p$ is finitely block (almost isometrically) represented in $X$, i.e. for any
$\e>0$ and any $n\in\N$ there is a normalized block sequence $x_1<\dots<x_n$ in $X$ which
is $(1+\e)$-equivalent to the unit vector basis of $\ell_p^n$. The set of all such $p$'s
is called the Krivine set of a given space $X$. Recall that the Krivine set of a
$p-$space $T[(\mc{A}_n,\theta_n)_{n\in\N}]$ contains only $p$ provided the sequence $(c_n)_n$ is
decreasing (Prop. 1.6 \cite{m}), while there are $p-$spaces with $p>1$
and 1 in the Krivine set (Prop. 1.8 \cite{m}).

\bt\label{kriv} Let $X$ be a $p-$space, $1\leq p<\infty$. Then the Krivine set of any
block subspace of $X$ contains $p$. \et

\bp By Lemma 1.5 \cite{m} Krivine set of $X$ is contained in $[1,p]$. We will show that
for any $N$ in any block subspace $Y$ of $X$ there is a normalized block sequence
$y_1,\dots, y_N$ with $\nrm{\sum_{n\in J}y_n}\leq D(\# J)^{1/p}$ for any $J\subset
\{1,\dots, N\}$ and for a universal constant $D$. Then by Corollary 5 \cite{mr} some $p'\geq
p$ is in the Krivine set of $Y$, thus $p'=p$ which ends the proof.

Fix $N\in \N$ and $r\leq p$ in the Krivine set of a block subspace $Y$. Take a normalized
block sequence $(y_n)_1^N$ and sequence $(m_n)_1^N\subset\N$ with

\bnum

\item $y_n$ is a $2-\ell_r-$average of length greater than $(2m_n)^r$, for any $n$,

\item $N\theta_{m_n}\leq 1/2^{n+2}$, for any $n$,

\item $\theta_{m_n}\sum_{i<n}\# \supp y_i\leq 1/2^{n+2}$ for any $n$,

\enum

We claim that $\nrm{\sum_{n=1}^Ny_n}\leq 99N^{1/p}$.

Take a functional $\phi$ with the tree analysis $(\phi_t)_{t\in T}$. By Lemma \ref{p-xi}
we can assume that this tree-analysis is comparable with $(y_n)$. Let $y=\sum_ny_n$.

Let $A=\{t\in T:\ \phi_t$ covers some $y_n$\}. Given any $t\in A$ let $I_t=\{
n=1,\dots,N:\ \phi_t$ covers $y_n\}$. Let $\theta_m$ be the weight of $\phi_t$.

Now consider cases:

\

CASE 1. $m\leq m_n$ for all $n\in I_t$.

For any $n\in I_t$ let $J_n=\{i\in\suc t:\ \supp \phi_i\subset \ran y_n\}$. By Lemma
\ref{p-av} we have $\sum_{i\in J_n}\phi_i(y)=\sum_{i\in J_n}\phi_i(y_n)\leq 8 (\#
J_n)^{1/q}$.

\

CASE 2. There is some $n\in I_t$ with $m>m_n$. Let $n_t$ be the
maximal element of $I_{t}$ with this property.

First let $L_t=\{n\not\in I_t:\ \supp y_n\cap \supp \phi\subset \supp \phi_t\}$. Notice
that for any $n\in L_t$ there is some $\phi_{t_n}$ - successor of $\phi_t$ so that $\supp
y_n\cap \supp \phi\subset \supp \phi_{t_n}$. Hence
$$
\phi_t(\sum_{n\in L_t}y_n)\leq \theta_{m_{n_t}}(\sum_{n\in L_t}\phi_{t_n}(y_n))\leq
N\theta_{m_{n_t}}\leq 1/2^{n_t+2}
$$
Hence $\phi(\sum_{t\in A,n\in L_t}y_n)\leq 1/4$, thus we can erase this part for all $t$
with error $1/4$.

Now notice that by condition 3. we have
$$
\phi_t(\sum_{n\in I_t, n<n_t}y_n)\leq \theta_{m_{n_t}}\sum_{n<n_t}\# \supp y_n \leq 1/
2^{n_t+2}
$$
so we can again erase this part for all $t$ with error $1/4$.

Now we compare the parts $\phi_t(y_{n_t})$ and $\phi_t(\sum_{n\in
I_t, n>n_t}y_n)$. If the first one is bigger, then set
$J_{n_t}=\{t\}$. Obviously we have $\phi_t(y)\leq (\#
J_{n_t})^{1/q}$. If the second one is bigger, then erase the first
one and proceed as in Case 1 getting $\sum_{i\in J_n}\phi_i(y)\leq
8 (\# J_n)^{1/q}$ for each $n\in I_t, n>n_t$.

Now by Lemma \ref{p-est} we have $\phi(\sum_{n=1}^N y_n)\leq
2\cdot 8 N^{1/p}+\frac{1}{2}$, therefore
$\nrm{\sum_{n=1}^Ny_n}\leq 6\cdot 16 \frac{1}{2}
N^{1/p}=99N^{1/p}$ which ends the proof of Theorem. \ep

Now we work under more restrictive assumptions on $(c_n)$, giving
additional information on representability of $\ell_p$ in
$p-$spaces.  Recall that a space with a basis is called
$\ell_p-$\textit{asymptotic} provided any normalized block
sequence $n\leq x_1<\dots<x_n$, $n\in\N$, is $C-$equivalent to the
unit vector basis of $\ell_p^n$, with some universal constant $C$.
Now we show a result generalizing Corollary 3.8 \cite{jko}.

\bpr Let $X$ be a $p-$space,  $1<p<\infty$, with $\inf_n c_n=c>0$. Then $X$ is saturated
with $\ell_p-$asymptotic subspaces. \epr

To prove the Proposition it is enough to show the following Lemma (well-known in case of
Tzafriri space) similar to Proposition 1.6 \cite{m}.

\bl\label{tz-est}

Assume $\inf_n c_n=c>0$. Take a normalized block sequence $y_1<\dots<y_N$ in $X$. Then
for any $(a_n)_{n=1}^N$ we have $\nrm{\sum_na_ny_n}\leq 6c^{-1}(\sum_na_n^p)^{1/p}$.

\el

\bp

Take any norming functional $\phi$ with the tree-analysis
$(\phi_t)_{t\in T}$, by Lemma \ref{p-xi} we can assume that the
tree-analysis is comparable with $(y_n)$. For any $n$ pick $t_n\in
T$ such that $\phi_{t_n}$ covers $y_n$ and let $(\phi_i)_{i\in
J_n}$ be the family of immediate successors of $\phi_{t_n}$ with
support contained in $\supp y_n$. By assumption and definition of
the norm we have $\sum_{i\in J_n}\phi_i(y_n)\leq c^{-1}(\#
J_n)^{1/q}$ for any $n$. Now we apply Lemma \ref{p-est} obtaining
$\phi(\sum_na_ny_n)\leq 6c^{-1}(\sum_na_n^p)^{1/p}$. \ep

By this Lemma and the fact that for any normalized $y_1<\dots<y_N$ also
$\nrm{\sum_ny_n}\geq cN^{1/p}$, applying Theorem 3.7 \cite{jko}, we obtain saturation of $X$ by $\ell_p-$asymptotic subspaces.

Recall that by proof of Theorem 2.1 \cite{m}, if $c_n\searrow 0$,
then $X$ does not contain a $\ell_p-$asymptotic subspace.

The next two observations should be compared with Proposition 3.9
\cite{ot} and Example 5.12 \cite{otw}, the idea of the proofs is
analogous.

\bpr Let $X$ be a $p-$space, $1<p<\infty$, with $\prod_n c_n=d>0$
or $c_n\nearrow 1$. Then $X$ is saturated with $\ell_p$. \epr

\bp By quasiminimality (Theorem \ref{p-quasi}) it is enough to show that there is a
subspace isomorphic to $\ell_p$. Since we have for any normalized block sequence the
upper $\ell_p-$estimate by Lemma \ref{tz-est}, it is enough to produce a normalized block
sequence with lower $\ell_p-$estimate. We present the proof in case $\prod_n c_n=d>0$,
the reasoning can be easily adapted to the case $c_n \nearrow 1$.

Take increasing $(i_n)_n\subset\N$ with $i_{n+1}>5i_n$ and $\sum_n2^{-i_n}\leq d/4$, and
take the partition $I_1<I_2<\dots$ of $\N$ with $\# I_n=2^{i_n}$. Notice that by Remark
\ref{rem-pnorm} we have $\nrm{\sum_{i\in I_n}e_i}\leq 2^{i_n/p}$.

Let $x_n=(\nrm{\sum_{i\in I_n}e_i})^{-1}\sum_{i\in I_n}e_i$ for any $n$. We claim that
$(x_n)$ satisfies the lower $\ell_p-$estimate.

Take $(j_n)_{n=1}^N\subset \N$ with $j_n\leq i_n$ and consider the vector
$x=\sum_n2^{-j_n/p}x_n$. Let $k_n=j_n+i_n$ for any $n$, then $k_{n+1}>k_n$. Then the
formula
$$
\phi=\sum_{n=1}^N\sum_{i\in I_n}(c_{2^{k_1}}c_{2^{k_2-k_1}}\dots
c_{2^{k_n-k_{n-1}}})2^{-k_n/q}e_i
$$
defines a functional from the norming set with a tree-analysis all of whose elements have
weights in the set $\{\theta_{2^{k_1}}, \theta_{2^{k_2-k_1}},\dots ,
\theta_{2^{k_n-k_{n-1}}}:\ n=1,\dots, N\}$. Moreover $k_{n+1}-k_n>k_n-k_{n-1}$ for any
$n$, hence
$$
\nrm{x}\geq\phi(x)\geq \sum_n\left(\prod_kc_k\right)\left( 2^{-k_n/q}2^{-j_n/p}\# I_n/\nrm{\sum_{i\in
I_n}e_i}\right)\geq d\sum_n 2^{-j_n}
$$
For any scalars $(a_n)_{n=1}^N\subset [0,1]$ with $\sum_n a_n^p=1$ pick $(j_n)$ with
$2^{(-j_n-1)/p}\leq a_n\leq 2^{-j_n/p}$. Let $J=\{ n\in\N: j_n\geq i_n\}$. Then by the
choice of $i_n$ we have
$$
\sum_{n\in J}a_n^p\leq\sum_{n\in J}2^{-j_n}\leq \sum_n 2^{-i_n}\leq d/4
$$
By the previous argument we have
$$
\nrm{\sum_na_nx_n}\geq \nrm{\sum_{n\not\in J}a_nx_n}\geq d\sum_{n\not\in J}2^{-j_n-1}\geq
d\sum_na_n^p/2 - d/4=d/4
$$
Hence we proved the lower $\ell_p-$estimate for $(x_n)$ and thus
the proposition. \ep

\bpr Let $X$ be a $p-$space, $1\leq p<\infty$, with $\sup_{n}c_{n}=d<1$. Then $X$ does not contain $\ell_p$. \epr

\bp Given a normalized block basis $(x_j)$ and $k\in\N$ define a new norm
$\nrm{\cdot}_{k,(x_j)}$ on $[x_j]$ in the following way: for any $x\in [x_j]$ let
$\nrm{x}_{k,(x_j)}$ be the supremum of $\phi(x)$ over all norming functionals $\phi$ with
a tree-analysis $(\phi_t)_{t\in T}$ such that for any $t\in T$ of level $\leq k$ and any
$j$ either $\supp\phi_t\cap \supp x_j=\emptyset$ or $\ran \phi_t\supset\supp x_j$.

\bl\label{not-p} Let $(x_j)$ be a normalized block basis in $X$, fix $\e>0$ and $k\in
\N$. Then there is $x\in [x_j]$ with $\nrm{x}=1$ such that $\nrm{x}_{k,(x_j)}>1-\e$. \el

\bp[Proof of Lemma] We will prove it first for $k=1$. Take a normalized block sequence
$(y^1_n)_{n=1}^N$ of $(x_j)$ such that

\bnum

\item $y^1_n=\sum_{j\in F_n}x_j/\nrm{\sum_{j\in F_n}x_j}$, where $\# F_n=m_n$ for any
$n$,

\item $m_n\sum_{i>n}1/\theta_{m_i}m_i\leq \e/4$ for any $n$,

\item $\theta_{m_n}\sum_{i<n}\# \supp y_i^1\leq \e/4$ for any $n$,

\item $1/N\theta_N\leq \e/8$.

\enum

Let $y^2=\sum_{n=1}^Ny_n^1/\nrm{\sum_{n=1}^Ny^1_n}$. Take $(E_i)_{i=1}^l$ with
$\nrm{y^2}=\theta_l\sum_i\nrm{E_iy^2}$. We claim that with error at most $\e$ we can
assume that for any  $n$ there is at most one $i$ with $\supp x_n\cap E_i\neq \emptyset$.
Let $J_n=\{ j\in F_n:\ \supp x_j$ intersects at least two $E_i\}$.

Assume that there is some $n_0$ such that $m_{n_0}<l\leq m_{n_0+1}$. Then we let $w_1$ be
$y^2$ restricted to supports of $y^1_n$ for $n<n_0$,  $w_2$ be $y^2$ restricted to
supports of $y^1_{n_0}$ and $y^1_{n_0+1}$. Then $\theta_l\sum_i\nrm{E_iw_1}\leq \e/4$ by
condition 3., $\theta_l\sum_i\nrm{E_iw_2}\leq 2/N\theta_N\leq \e/4$ by condition 4., and
by condition 1. we have
$$
\sum_{n>n_0+1}\sum_{j\in J_n}\nrm{y_{|\supp x_j}}\leq
m_{n_0+1}\sum_{n>n_0+1}1/\theta_{m_n}m_n\leq \e/4
$$
Any other case is a simpler modification of the reasoning above.

For any $k>1$ we iterate the procedure above, by taking as
$(y^j_n)$ suitable averages of $(y_n^{j-1})$ for some $\e_j$,
$j=2,\dots,k$, with $\prod_{i=1}^k(1-\e_i)>1-\e$. \ep

Now we go back to the proof of the Proposition. Assume there is a normalized block basis
$(x_j)$ $D$-equivalent to the unit vector basis of $\ell_p$. Take $k\in\N$ and by the
Lemma \ref{not-p} pick a normalized $x=\sum_ja_jx_j$ with $\nrm{x}_{k,(x_j)}>1/2$. Let
$$
\nrm{x}_{k,(x_j)}=\sum_i \frac{c_{n_1^i}c_{n_2^i}\dots
c_{n_{k_i}^i}}{(n_1^i\cdot\dots\cdot n_{k_i}^i)^{1/q}} \nrm{E_ix}
$$
for some $E_i$ such that the sets $J_i=\{j:\ \supp x_j\cap E_i\neq\emptyset\}$ are
pairwise disjoint and some $k_i\leq k$ and $(n^i_1,n^i_2,\dots n^i_{k_i})_i$ with $\sum_i 1/(n_1^i\cdot\dots\cdot
n_{k_i}^i)\leq 1$. Then by H\"older inequality
$$
1/2\leq d^k\sum_i 1/(n_1^i\dots n_{k_i}^i)^{1/q} D(\sum_{j\in
J_i}a_j^p)^{1/p}\leq d^kD(\sum_ja_j^p)^{1/p}\leq d^{k}D^2
$$
which for large $k$ gives contradiction. \ep

\section{Quasiminimality in $p-$spaces}

In this section we prove the following

\bt\label{p-quasi} Any $p-$space $X$, $1\leq p<\infty$, is quasiminimal. \et

We show also a stronger quasiminimality property of $X$ and its dual space under
additional assumption on $(c_n)$. The equivalent sequences used in the proof of Theorem
\ref{p-quasi} consist of averages of $\ell_p-$averages, whose properties are examined in
next Lemmas.

\bl \label{p-eq} Fix $N\in \N$. Take finite normalized block
sequences $(y_n)_1^N$ and $(z_n)_1^N$ in $X$ and sequence
$(m_n)_1^N\subset\N$ satisfying the following conditions:

\bnum

\item each $y_n$ and $z_n$ are $2-\ell_p-$averages of length $N_n\geq (2m_n)^p$,

\item $N\theta_{m_n}\leq 1/2^{n+5}$ for any $n$,

\item $\theta_{m_n}\sum_{i<n}\# \supp y_i\leq 1/2^{n+5}$, and $\theta_{m_n}\sum_{i<n}\#
\supp z_i\leq 1/2^{n+5}$ for any $n$,

\enum

Then $(y_n)_1^N$ and $(z_n)_1^N$ are $C$-equivalent, for some universal constant C.
\el

\bp We show that $(z_n)_1^N$ dominates $(y_n)_1^N$. Fix scalars $0\leq a_1,\dots,a_N\leq
1$ so that the vector $y=\sum_{n=1}^Na_ny_n$ has norm 1, and put $z=\sum_{n=1}^Na_nz_n$.
The proof resembles the reasoning in the proof of Theorem \ref{kriv}.

By Lemma \ref{p-xi} we can take a norming functional $\phi$ with a tree analysis $(\phi_t)_{t\in T}$ comparable with $(y_n)$
and $\phi(y)\geq 1/6$. We will pick $\psi$ so that $\nrm{\psi}\leq 1$
and $\psi (z)\geq 2^{-11}$.

Let $A=\{t\in T:\ \phi_t$ covers some $y_n$\}. Given any $t\in A$ let $I_t=\{
n=1,\dots,N:\ \phi_t$ covers $y_n\}$. Let $\theta_m$ be the weight of $\phi_t$.

Now consider cases for each $t\in A$:

\

CASE 1. $m\leq m_n$ for all $n\in I_t$. For any $n\in I_t$ let $J_n=\{i\in\suc t:\ \supp
\phi_i\subset \ran y_n\}$. By condition 1. and Lemma \ref{p-av} pick functionals
$(\psi_i)_{i\in K_n}$, with $\# K_n\leq \# J_n$ so that $\sum_{i\in K_n}\psi_i(z_n)\geq
\sum_{i\in J_n}\phi_i(y_n)2^{-6}$ and change the tree-analysis of $\phi$ by putting
$(\psi_i)_{i\in K_n}$ instead of $(\phi_i)_{i\in J_n}$ for any $i=1,\dots,k_n$. Then we
have $2^6\psi_t(\sum_{n\in I_t}a_nz_n)\geq \phi_t(\sum_{n\in I_t}a_ny_n)$.

\

CASE 2. there is some $n\in I_t$ with $m>m_n$. Let $n_t$ be maximal with this property.

First let $L_t=\{n\not\in I_t:\ \supp y_n\cap \supp \phi\subset \supp \phi_t\}$. Notice
that for any $n\in L_t$ there is some $\phi_{t_n}$ - a successor of $\phi_t$ so that
$\supp y_n\cap \supp \phi\subset \supp \phi_{t_n}$. Hence
$$
\phi_t(\sum_{n\in L_t}a_ny_n)\leq \theta_{m_{n_t}}(\sum_{n\in L_t}\phi_{t_n}(y_n))\leq
N\theta_{m_{n_t}}\leq 1/2^{n_t+5}
$$
Hence $\phi(\sum_{n\in L_t,t\in B}a_ny_n)\leq 1/2^5$, where $B=\{t\in A: m>m_k$ for some
$k\in I_t\}$, thus we can erase this part with error $1/2^5$.

Now notice that by condition 2. we have
$$
\phi_t(\sum_{n\in I_t, n<n_t}a_ny_n)\leq \theta_{m_{n_t}}\sum_{n<n_t}\# \supp y_n \leq 1/
2^{n_t+5}
$$
so we can again erase this part for all $t$ with error $1/2^5$.

Now we compare the part $\phi_t(a_{n_t}y_{n_t})$ and
$\phi_t(\sum_{n\in I_t, n>n_t}a_ny_n)$. If the first one is
bigger, then replace $\phi_t$ in the tree-analysis of $\phi$ by
$\psi_t$ satisfying $\psi_t (z_{n_t})=1$. If the second one is
bigger, then erase the first one and proceed as in Case 1.

Altogether after all this changes we obtained a new norming functional $\psi$ with
$\psi(\sum_na_nz_n)\geq 2^{-11}$. Indeed, from the above inequalities we have
\begin{equation}
 \frac{1}{6}\leq \phi(\sum_{n}y_{n})\leq \frac{2}{2^5}+2\cdot
 2^{6}\psi(\sum_{n}a_{n}y_{n})\
\Rightarrow\ \psi(\sum_{n}a_{n}y_{n})\geq (\frac{1}{6}-\frac{1}{2^4})\frac{1}{2^7}.
\end{equation}
Therefore $C= 2^{11}$ satisfies the assertion of the Lemma. \ep

\bl \label{p-eq2} Fix $N\in\N$. Take $(y_n)_1^N$, $(z_n)_1^N$ as in Lemma \ref{p-eq} and
let $y=y_1+\dots+y_N$, $z=z_1+\dots+z_N$. Then for $m\in\N$, any $E_1<\dots<E_m$
there are $F_1<\dots<F_l$ with $l\leq m$
\begin{equation}\label{pick}
\theta_m(\nrm{E_1y}+\dots+\nrm{E_my})\leq C'\theta_m(\nrm{F_1z}+\dots+\nrm{F_lz})+6
\end{equation}
for some universal constant $C'$.
\el

\bp We can assume (up to multiplying by 3) that the sequence $(E_j)$ are comparable with
the sequence $(y_n)$. Let $J=\{ j=1,\dots,m:\ \ran E_j\supset\supp y_n$ for some $n\}$
and for $j\in J$ let $I_j=\{n:\ \supp y_n\cap\bigcup_iE_i\subset E_j\}$. Then for any
$j\in J$ we put $F_j=\bigcup_{n\in I_j}\supp z_n$ and we have by Lemma \ref{p-eq}
$$
\sum_{j\in J}\nrm{E_j\sum_{n\in I_j}y_n}\leq C\sum_{j\in J}\nrm{F_j\sum_{n\in I_j}z_n}
$$
Symmetrically let $L=\{n:\ \ran y_n\supset E_j$ for some $j\}$ and
for any $n\in L$ let $J_n=\{j:\ E_j\subset \ran y_n\}$. Consider
cases:

\

CASE 1. $m\leq m_n$ for any $n$. Then by the choice of $(y_n),(z_n)$ and Lemma \ref{p-av}
there are some $(F_j)_{j\in K_n, n\in I}$, $\# K_n\leq \# J_n$ with
$$
\sum_{n\in L}\sum_{j\in J_n}\nrm{E_jy_n}\leq 2^6 \sum_{n\in L}\sum_{j\in K_n}\nrm{F_jz_n}
$$
and $\{F_j\}_{j\in J}\cup\{F_j\}_{j\in K_n, n\in L}$ are the
desired sets.

\

CASE 2. For some $n$ we have $m>m_n$. Let $n_0$ be maximal with this property. Then
$$
\theta_m\sum_{n\in L,n<n_0}\sum_{j\in J_n}\nrm{E_jy_n}\leq \theta_{m_{n_0}}\sum_{n<n_0}\#
\supp y_n\leq 1/2^{n_0}.
$$
Notice that
$$
\theta_m\sum_{j\in J_{n_0}}\nrm{E_jy_{n_0}}\leq \nrm{y_{n_0}}\leq 1
$$
For $(E_j)_{j\in J_n,n\in L,n>n_t}$ choose $(F_j)$ as in the Case
1. Thus we prove the Lemma with constant $C'=3\cdot 2^6C$. \ep

\bp[Proof of the Theorem \ref{p-quasi}] Take any block subspaces $Y,Z$ of $X$ and pick by
Theorem \ref{kriv} $(y_j)\subset Y$, $(z_j)\subset Z$ so that each
$y_j=(y^j_1+\dots+y^j_{N_j})/\nrm{y^j_1+\dots+y^j_{N_j}}$,
$z_j=(z^j_1+\dots+z^j_{N_j})/\nrm{z^j_1+\dots+z^j_{N_j}}$, were $(y_n^j)$, $(z_n^j)$ are
as in Lemmas \ref{p-eq}, \ref{p-eq2} and $N_j$ are so big that
$\nrm{y^j_1+\dots+y_{N_j}^j}, \nrm{z^j_1+\dots+z_{N_j}^j}\geq 2^{j+7}$ (for example
$N_j\theta_{N_j}\geq 2^{j+7}$). Take any scalars $(a_j)\subset [0,1]$ with
$\nrm{\sum_ja_jy_j}=1$. By Lemma \ref{p-xi} we pick a norming functional $\phi$ with
$\phi(\sum_ja_jy_j)\geq 1/6 $ and a tree-analysis $(\phi_t)_{t\in T}$ comparable with
$(y_n)$. We will pick norming $\psi$ with $\psi(\sum_j a_jz_j)\geq 1/2^4C'$.

 For any $j$ pick $t_j$ so that $\phi_{t_j}$ covers $y_j$. For
any $j$ let $\phi_1^j<\dots<\phi^j_{k_j}$ be all the successors of $\phi_{t_j}$ whose
supports are contained in $\supp y_j$. Then by Lemma \ref{p-eq2} pick for any $j$ suitable
$F_1^j<\dots<F_{l_j}^j$ with $l_j\leq k_j$ satisfying \eqref{pick}. Replace in the tree-analysis of
$\phi$ each $\phi_i^j$ by $\psi_i^j$ so that $\nrm{F_i^jz_j}=\psi_i^j(z_j)$ and $\supp
\psi_i^j\subset \supp z_j$. Then for any $j$
$$
\theta_{k_j}\sum_{i=1}^{k_j}\phi^j_i(y_j)\leq C'\theta_{k_j}\sum_{i=1}^{l_j}\psi_i^j(z_j)+6/2^{j+7}
$$
hence
$$
1/6\leq \phi(\sum_ja_jy_j)\leq C'\psi(\sum_ja_jz_j)+1/2^4
$$
which ends the proof, since $\psi$ is from the norming set. \ep

\br Notice that the proof provides in every block subspace of $X$ some subsequentially
minimal block subspace. Indeed any block sequence of arbitrary long averages of vectors
of the form $(y_1+\dots+y_N)/\nrm{y_1+\dots+y_N}$, where $y_1,\dots,y_N$ are as in Lemma
\ref{p-eq} for arbitrary big $N$ and $m_1,\dots,m_N$, spans a subsequentially minimal
subspace.

Moreover by  Lemma \ref{p-eq} sequences of long $\ell_p-$averages generates equivalent
spreading models. In case $p=1$ the proof can be easily adapted to show that the
spreading models of sequences of long $\ell_1-$averages are equivalent to the unit vector
basis of the considered space $X$.

\er

Now we consider special case: $\inf_n c_n=c>0$. An example of such
space is Tzafriri space $T[(\mc{A}_n,c/\sqrt{n})_n]$. Recall that
a Banach space $X$ with a basis belongs to Class 2 defined by
T.Schlumprecht \cite{s}, if any block subspace $Y$ of $X$ contains
normalized block sequences $(x_n)$, $(y_n)$ such that the mapping
carrying each $x_n$ to $y_n$ extends to a bounded strictly
singular operator, i.e. operator whose no restriction to an
infinite dimensional subspace is an isomorphism.

\bt\label{tz-quasi} Let $X$ be a $p-$space, $1<p<\infty$, with $\inf_n c_n=c>0$. Then

\bnum

\item the dual space $X^*$ is quasiminimal,

\item the spaces $X$ and $X^*$ do not contain a subspace of Class 2.

\enum

\et

\br By the above theorem the class of $p-$spaces with $\inf_nc_n>0$ might provide an
example of a Banach space containing neither a subspace of Class 2 nor a subspace of
Class 1 (recall that a space with a basis is of Class 1 if any its normalized block
sequence admits a subsequence equivalent to some subsequence of the basis, \cite{s}),
answering the question (Q4), \cite{s}.

\er

In order to prove the theorem we need the following

\bl\label{tz-eq}  Assume $\inf_n c_n=c>0$. Take normalized block sequences $(y_n)$ and
$(z_n)$ in $X$ such that $\nrm{z_n}_\infty\leq c/16(\# \supp y_n)^2$ for any $n$. Then
$(z_n)$ $D-$dominates $(y_n)$ with some universal constant $D$.

\el

\bp[Proof of Lemma \ref{tz-eq}] We prove two Claims:

\bcl\label{partition} Fix $m\in \N$ and $1>\de>0$. Then for any $z$ with $\nrm{z}\geq
1/2$, $\nrm{z}_{\infty}<\de/8m^2$ there are $F_1<\dots<F_m$ such that
$$
1-\de\leq\nrm{F_iz}/\nrm{F_jz}\leq 1+\de, \ \ \ i,j=1,\dots,m
$$

\ecl

\bp[Proof of Claim \ref{partition}] Fix $\e>0$. We prove induction on $m$ that for any
$z$ with $\nrm{z}_{\infty}<\e$ and $m\leq \# \supp z$ there are $F_1<\dots<F_m$ such that
$\sup_i\nrm{F_iz}\leq \inf_i\nrm{F_iz}+2m\e$ and $F_1\cup\dots\cup F_m=\supp z$.

For $m=2$ the claim easily follows. Assume the claim holds for $m-1$. Take $z$ with
$\nrm{z}_\infty<\e$ and $\#\supp z\geq m$. We can assume that $\supp z=\{1,\dots,J\}$ for
some $J\geq m$ to simplify the notation. For any $m-1\leq j\leq J$ let $z_j$ be the
restriction of $z$ to the interval $\{1,\dots,j\}$. By the inductive assumption for $m-1$
pick suitable $F_1^j<\dots<F_{m-1}^j$ for vectors $z_j$, $m-1\leq j\leq J$.

Notice that $\sup_i\nrm{F_i^jz_j}\leq \sup_i\nrm{F_i^{j-1}z_{j-1}}+(2m-1)\e$. Indeed, if
not, then by the choice of $F_i^j$'s we have $\nrm{F_i^jz_j}>\nrm{F_k^{j-1}z_{j-1}}+\e$
for any $1\leq i,k\leq m-1$, hence $\max F_i^j>\max F_i^{j-1}+1$ for any $i\leq m-1$,
which by assumption on $\nrm{z}_\infty$ contradicts $\max\supp z_j=\max\supp z_{j-1}+1$.

Notice also that $\nrm{z-z_j}+\e\geq \nrm{z-z_{j-1}}$ for any $m-1\leq j\leq J$.

Now if $\nrm{z-z_{m-1}}<\e$, then $\{ 1\},\dots, \{ m-1\}, \supp (z-z_{m-1})$ satisfy the
claim.

Otherwise consider function $\xi(j)=\nrm{z-z_j}-\sup_i\nrm{F_i^jz_j}$ for $m-1\leq j\leq
J$. We have that
$$
\xi(m-1)=\nrm{z-z_{m-1}}-\sup_{i}\nrm{F_{i}^{m-1}z_{m-1}}\geq
\nrm{z-z_{m-1}}-\nrm{z}_{\infty}>0
$$
and $\xi(J)< 0$ and by previous observations
$\xi(j)+2m\e\geq \xi(j-1)$ for any $m-1<j\leq J$. Take the minimal $J\geq j_0> m-1$ with
$\xi(j_0)\leq 0$, then $0<\xi(j_0-1)\leq \xi(j_0)+2m\e$, hence
$$
\sup_i\nrm{F_i^{j_0}z_{j_0}}-2m\e\leq\nrm{z-z_{j_0}}\leq \sup_i\nrm{F_i^{j_0}z_{j_0}}
$$
thus if we put $F_i=F_i^{j_0}$ for any $1\leq i<m$ and $F_m=\supp (z-z_{j_0})$, then we
have by the inductive assumption and choice of $j_0$ that $\sup_{1\leq i\leq
m}\nrm{F_iz}\leq \inf_{1\leq i\leq m}\nrm{F_iz}+2m\e$.

\

Take $(F_i)_1^m$ as above. Take $z$ as in Claim, then also
$\#\supp z\geq m$. Notice that $\sup\nrm{F_iz}\geq 1/2m$. Hence
$\inf\nrm{F_iz}/\sup\nrm{F_iz}\geq 1 - 4m^2\e$ thus also
$\sup\nrm{F_i z}/\inf\nrm{F_iz}\leq 1/(1-4m^2\e)$ which ends the
proof. \ep

\bcl\label{tz-eq2} Assume $\inf_nc_n=c>0$. Take any normalized
$y,z\in X$, assume that $\nrm{z}_\infty\leq c/16(\# \supp y)^2$.
Then for $m\in\N$, any $E_1<\dots<E_m$ there are $F_1<\dots<F_m$
with
$$
\nrm{E_1y}+\dots+\nrm{E_my}\leq 9c^{-2}(\nrm{F_1z}+\dots+\nrm{F_mz})
$$
\ecl

\bp[Proof of Claim \ref{tz-eq2}] Notice that by definition of the norm we have
$\nrm{E_1y}+\dots+\nrm{E_my}\leq m^{1/q}/c$. Of course we can assume that $m\leq \# \supp
y$. Now we pick $F_1<\dots<F_m$ by the Claim \ref{partition} for $\de=1/2$. Then by Lemma
\ref{tz-est} we have
$$
c\nrm{z}/6\leq (\sum_i\nrm{F_iz}^p)^{1/p}\leq
m^{1/p}3\inf_i\nrm{F_iz}/2
$$
therefore $\nrm{F_iz}\geq c/9m^{1/p}$ for any $i$, hence
$$
\sum_i\nrm{F_iz}\geq m^{1/q}c/9
$$
which end the proof of the Claim. \ep

This Claim shows in fact that Lemma \ref{p-eq2} holds in $X$ without error for much
larger class of vectors provided the special assumption on $(c_n)$ holds. Hence coming
back to the proof of Lemma \ref{tz-eq} we proceed as in the proof of Theorem
\ref{p-quasi}, the situation is now simpler than in the general case, enabling us to pass
to the dual space. Take any $a_1,\dots,a_N\in [0,1]$. Take any functional $\phi$ from the
norming set with a tree-analysis $(\phi_t)_{t\in T}$, assume by Lemma \ref{p-xi} that
$(\phi_t)_{t\in T}$ is comparable with $y_1,\dots,y_N$. We will show that there is a
functional $\psi$ from the norming set such that $\phi(\sum_na_ny_n)\leq 9c^{-4}
\psi(\sum_na_nz_n)$. For any $n$ take $t_n$ such that $\phi_{t_n}$ covers $y_n$ and let
$(\phi_i)_{i\in J_n}$ be immediate successors of $\phi_{t_n}$ with supports contained in
$\supp y_n$. Using Claim \ref{tz-eq2} find $(\psi_i)_{i\in J_n}$ with
$$
\sum_{i\in J_n}\phi_i(y_n)\leq 9c^{-2}\sum_{i\in J_n}\psi_i(z_n)
$$
and replace in the tree-analysis of $\phi$ all $(\phi_i)_{i\in
J_n}$ with $(\psi_i)_{i\in J_n}$ obtaining in such a way a norming
functional $\psi$ with the desired property. \ep

\bp[Proof of Theorem \ref{tz-quasi}]

We proceed to the proof of the first part. We need the following claims:

\bcl\label{al-de-vep} For any $\de,\al>0$ there is some $\e>0$ such that for any $\phi$
with $\nrm{\phi}_{\infty}<\e$ and any $x\in X$ with $\nrm{x}\leq 1$ and $\phi(x)>\al$
there is some $I\subset\N$ so that $\abs{\phi(x_I)}\geq\abs{\phi(x)}/2$ and
$\nrm{x_I}_{\infty}<\de$.

\ecl

\bp[Proof of Claim \ref{al-de-vep}]

Fix $\de,\al>0$. Take $\phi$ with $\nrm{\phi}_{\infty}<\e$ for some $\e>0$ and $x$ with
$\nrm{x}\leq 1$ and put $J=\{i\in \supp x :\ \abs{x(i)}>\de\}$. Then $c(\#
J)^{1/p}\de\leq \nrm{x}\leq 1$ from the definition of the norm in $X$, hence
$$
\abs{\phi(x_J)}=\sum_{j\in J}\abs{\phi(j)x(j)}<\e\# J\leq \e /(c\de)^p
$$
Thus $\e=(c\de)^p\al/2$ and $I=\N\setminus J$ satisfy the desired property. \ep

\bcl\label{domin-dual} Take normalized block sequences $(\phi_n)$, $(\psi_n)$ in the
norming set of $X$ such that $\nrm{\psi_n}_{\infty}<\e_n$, where each $\e_n$, $n\in\N$,
is chosen by Claim \ref{al-de-vep} for $\de=c/4^{n+3} (\# \supp \phi_n)^2$ and
$\al=1/4^n$.

Then $(\phi_n)$ $4D$-dominates $(\psi_n)$, with $D$ as in Lemma \ref{tz-eq}.

\ecl

\bp[Proof of Claim \ref{domin-dual}]

Take $a_1,\dots,a_N$ and consider $\phi=a_1\phi_1+\dots+a_N\phi_N$,
$\psi=a_1\psi_1+\dots+a_N\psi_N$. Let $\nrm{\psi}=1$ and take $z\in X$ with $\nrm{z}=1$
and $\psi(z)=1$. We can assume, by unconditionality of the unit vector basis, that all
coefficients of $\phi$ and $z$ are non-negative, as well as $a_1,\dots,a_N$.

Let $J=\{n:\ \psi_n(z)\leq 1/4^n\}$. Then $\sum_{n\in J}a_n\psi_n(z)\leq 1/2$.

By the choice of $(\e_n)$ we can take for any $n\not\in J$ some
vector $w_n$ - a restriction of $z_{\supp\psi_n}$ - such that
$\nrm{w_n}_\infty\leq c/4^{n+3} (\# \supp \phi_n)^2$ and
$\nrm{w_n}\geq \psi_n(w_n)\geq \psi_n(z)/2\geq 1/4^{n+1}$.

Take $v_n=w_n/\nrm{w_n}$, $n\not\in J$, and notice that $\nrm{v_n}_\infty\leq c/16(\#
\supp \phi_n)^2$. Now pick normalized block sequence $(y_n)$ with $\phi_n(y_n)=1$, $\supp
y_n\subset \supp \phi_n$. By Lemma \ref{tz-eq} we have that $(v_n)$ $D$-dominates
$(y_n)$, hence the vector $y=\sum_{n\not\in J}\nrm{w_n}y_n$ satisfies $\nrm{y}\leq
D\nrm{\sum_{n\not\in J} w_n}\leq D\nrm{z}=D$. Therefore
$$
D\nrm{\phi}\geq \phi(y)=\sum_{n\not\in J} a_n\nrm{w_n}\geq\sum_{n\not\in J}
a_n\psi_n(w_n)\geq\sum_{n\not\in J} a_n\psi_n(z)/2\geq 1/4
$$ \ep

Now take any block subspace $Y,Z\subset X^*$. Pick inductively infinite normalized block
sequences $(\phi_n)\subset Y$, $(\psi_n)\subset Z$ such that

\bnum

\item $\nrm{\psi_n}_\infty<\e_n$, where $\e_n$ is chosen by Claim \ref{al-de-vep}
for $\de=c/4^{n+3} (\# \supp \phi_n)^2$ and $\al=1/4^n$, for any $n$,

\item $\nrm{\phi_{n+1}}_{\infty}<\e'_n$, where $\e'_n$ is chosen by Claim
\ref{al-de-vep} for $\de=c/4^{n+3} (\# \supp \psi_n)^2$ and $\al=1/4^n$, for any $n$.

\enum

Then by Claim \ref{domin-dual}, $(\phi_n)$ $4D$-dominates
$(\psi_n)$ and $(\psi_n)$ $4D$- dominates $(\phi_{n+1})$. Now
split $\N$ into intervals $I_1<I_2<\dots$ such that
$\nrm{\sum_{n\in I_j}\phi_n}\geq 2^j$ and $\nrm{\sum_{n\in
I_j}\psi_n}\geq 2^j$ (this is possible, since $X^*$ does not
contain $c_0$). The following Lemma ends the proof of the first
part of the Theorem:

\bl \label{averages} Let $X$ be a Banach space. Take normalized 1-unconditional basic
sequences $(u_n), (v_n)\subset X$ so that $(u_n)$ dominates $(v_n)$, $(v_{n+1})$
dominates $(u_n)$ and sequence $I_1<I_2<...$ of intervals
of $\N$ such that
$\nrm{\sum_{n\in I_j}u_n}\geq 2^j$, $\nrm{\sum_{n\in I_j}v_n}\geq 2^j$. Let
$y_j=\sum_{n\in I_j}u_n/\nrm{\sum_{n\in I_j}u_n}$, $z_j=\sum_{n\in
I_j}v_n/\nrm{\sum_{n\in I_j}v_n}$ for any $j\in\N$. Then $(y_j)$ and $(z_j)$ are
equivalent.
\el
\bp[Proof of Lemma]
Let $D\geq 1$ be the domination constant. Notice first that
$$
\nrm{\sum_{n\in I_j}v_n}\leq D\nrm{\sum_{n\in I_j}u_n}$$
and
$$
\nrm{\sum_{n\in I_j}u_n}\leq D\nrm{\sum_{n\in
I_j}v_{n+1}}\leq D\nrm{\sum_{n\in I_j}v_n}+D\nrm{v_{\max I_j+1}}\leq 2D\nrm{\sum_{n\in I_j}v_n}
$$
Hence, by $(u_n)$ dominating $(v_n)$, for any scalars $(a_j)$ we
have
$$
\nrm{\sum_ja_jz_j}\leq 2D^2 \nrm{\sum_ja_jy_j}
$$
To see the other domination take $(a_j)$ with $\nrm{\sum_ja_jz_j}=1$. Then, by
$(v_{n+1})$ dominating $(u_n)$, we have
\begin{align*}
\nrm{\sum_ja_jy_j} &\leq \sum_j(|a_j|/\nrm{\sum_{n\in I_j}u_n}) +
\nrm{\sum_ja_j(\sum_{n\in I_j\setminus\{\max I_j\}}u_n)/\nrm{\sum_{n\in I_j}u_n}}
\\
&\leq 1+2 D^2\nrm{\sum_ja_jz_j}\leq 3 D^2
\end{align*}
\ep

In order to show the second part of the Theorem in case of space $X$ we prove the
following

\bcl Let $\inf_n c_n=c>0$. Take a normalized block sequences
$(u_n), (v_n)\subset X$ with $\nrm{u_n}_\infty\to 0$ and
$\nrm{v_n}_\infty\to 0$. Then if the mapping $v_n \to u_n$ extends
to a bounded linear operator $T$ then there is an infinite
dimensional subspace $Y$ of $[u_n]$ with $T|_Y$ an isomorphism.
\ecl
\bp[Proof of Claim] Take $(u_n)$, $(v_n)$ as above. Assume
that $T$ defined as above is strictly singular. Pick inductively
subsequence $(u_{i_n})$ and $(v_{j_n})$ such that
$\nrm{v_{j_n}}_\infty<c/16(\#\supp v_{i_n})$ and
$\nrm{u_{i_{n+1}}}_\infty<c/16(\#\supp v_{j_n})$. Then by Lemma
\ref{tz-eq} $(v_{j_n})$ dominates $(v_{i_n})$, hence also
$(u_{i_n})$, and $(u_{i_{n+1}})$ dominates $(v_{j_n})$. Denote the
operator carrying each $v_{i_n}$ to $u_{i_n}$ by $T'$. As a
composition of a strictly singular operator and a bounded
right-hand shift operator $T'$ is also strictly singular.

On the other hand we can consider averages $y_j=\sum_{n\in
I_j}u_{i_n}/\nrm{\sum_{n\in I_j}u_{i_n}}$, $z_j=\sum_{n\in
I_j}v_{i_n}/\nrm{\sum_{n\in I_j}v_{i_n}}$ for some sufficiently
long $(I_j)$. Then by Lemma \ref{averages} $(y_j)$, $(z_j)$ are
equivalent and moreover
$$
1/D'\leq \nrm{\sum_{n\in
I_j}u_{i_n}}/\nrm{\sum_{n\in I_j}v_{i_n}}\leq D'
$$
for some $D'\geq 1$. Hence $T'$ restricted to $[z_j]$ is an isomorphism, which gives
contradiction. \ep

The case of the dual space $X^*$ follows analogously, by using Claim \ref{domin-dual}.
\ep

\section{Quasiminimality in mixed Tsirelson spaces $T[(\mc{S}_n,\theta_n)_n]$}

Notice (\cite{ao}) that any mixed Tsirelson space $T[(\mc{S}_{k_n},\theta_n)_{n\in \N}]$,
with $\theta \to 0$, is isometric to $T[(\mc{S}_n,\widehat{\theta}_n)_{n\in\N}]$, where
$$
\widehat{\theta}_n=\sup\left\{\prod_{i=1}^l\theta_{n_i}:\ \sum_{i=1}^l n_i\geq n
\right\}, \ \ \ n\in N
$$
Notice that in such case the sequence $(\widehat{\theta}_n)$
satisfies $\widehat{\theta}_n\searrow 0$ and
$\widehat{\theta}_{n+m}\geq \widehat{\theta}_n\widehat{\theta}_m$
for any $n,m\in\N$. We will assume in the paper that we work in
such a setting.

\bd \cite{ao} A mixed Tsirelson space
$T[(\mc{S}_n,\theta_n)_{n\in\N}]$ is called \textit{regular} if
$\theta_n\searrow 0$ and $\theta_{n+m}\geq\theta_n\theta_m$ for
any $n,m\in\N$. \ed

\begin{notation}

Take a regular mixed Tsirelson space $T[(\mc{S}_n,\theta_n)_{n=1}^\infty]$. Then there
exists $\theta=\lim_n\theta_n^{1/n}=\sup_n\theta_n^{1/n}\in (0,1]$. We will use the
following notation: $\theta_n=c_n\theta^n$, with $c_n\in (0,1]$, $n\in\N$. If $c_n=1$ for
some $n\in\N$ then $T[(\mc{S}_n,\theta_n)]$ is isomorphic to Tsirelson space
$T[\mc{S}_1,\theta]$.

\end{notation}

\subsection{Two lemmas for the Schreier families}

We prove first technical lemmas needed in the sequel.

\bl\label{sch1}

Let $k,m\in\N$ and  $l\in\N$  such that  $km< 2^{l}$. Then
$$
(\mc{S}_{n}[\mc{A}_k])[\mc{A}_m] \subset \mc{A}_{l}[\mc{S}_{n}]\,\,\,\text{for any $n\in\N$}.
$$
\el \bp We prove the
result by induction on $n$. For $n=1$ we have
$$
[\mc{S}_1[\mc{A}_{k}]][\mc{A}_{m}] = \mc{S}_1[\mc{A}_{k}[\mc{A}_{m}]] = \mc{S}_1[\mc{A}_{km}]
$$
Let $j \leq m_1<\dots<m_j$ and
$$
m_i\leq \lambda_{(i-1)km+1} <  \dots < \lambda_{(i-1)km+km-1} < \lambda_{kmi} < m_{i+1}
$$
for every $i=1,\dots,j$. We set $G=\{\lambda_1,\dots, \lambda_{kmj}\}$. From the
assumptions we have that $j \leq G$.

For $i=1,\dots,l$ consider the sets
$$
F_{i}=\{2^{i-1}j,2^{i-1}j+1,\dots,2^{i}j-1\}\,\,\,\,
$$
Each of the sets $F_i$, $i\leq l$ belongs to the family $\mc{S}_1$ since $\# F_i \leq
\min F_i$. We may write $G\subset\bigcup_{i=1}^{l} G_i$ where $G_i =\{\lambda_{r-j+1} :\
r \in F_i\}$ for $i\leq l$. Then each $G_i$ is a spreading of $F_i$ and hence belongs to
the Schreier family $\mc{S}_1$. Thus $G\in \mc{A}_{l}[\mc{S}_1]$.

Assume that the result holds for $n$, i.e. $[\mc{S}_n[\mc{A}_{k}]][\mc{A}_{m}]\subset
\mc{A}_{l}[\mc{S}_{n}]$. Then we have that
$$
[\mc{S}_{n+1}[\mc{A}_{k}]][\mc{A}_{m}] = \mc{S}_1[\mc{S}_{n}[\mc{A}_{k}]][\mc{A}_{m}]
\subset \mc{S}_{1} [\mc{A}_{l}[\mc{S}_{n}]]\subset \mc{A}_{l}[\mc{S}_{n+1}]
$$
The last inclusion follows immediately from the associativity of
the operation $M[N]$ and the fact that
$\mc{A}_{l}[\mc{S}_n]=(\mc{S}_n)^{l}$. \ep

\bl\label{l3} Let $F\in\mc{S}_{m}$ and $(F_{i})_{i=1}^{d}$ be a
partition of $F$ into successive sets such that  for every  $i\leq
d$, $F_{i}\in\mc{S}_{m_i}\setminus\mc{S}_{m_{i}-1}$ for some
$m_{i}\leq m$. For every $i\leq d$ let $G_{i}\in\mc{S}_{m_i-1}$
(in case $m_{i}=0$ let $G_{i}\in \mc{S}_{0}$) with
$G_1<\dots<G_d$. Then we have the following:

\bnum

\item if $3\max F_{i}<\min G_{i+1}$ for every $i=1,\dots,d-1$,
then $\bigcup_{i=2}^{d}G_{i}\in \mc{S}_{m}$  and
$\bigcup_{i=1}^{d}G_{i}\in\mc{A}_{2}[\mc{S}_{m}]$.

\item if $2\max F_{i}<\min G_{i}$ for every $i=1,\dots, d$, then
$\bigcup_{i\leq d}G_{i}\in\mc{S}_{m}$. \enum
\el

\bp

We prove the first part by induction on $m$. For $m=1$ the
result is clear since $G_{i}\in\mc{S}_{0}$ and also for $i\geq 2$,
$G_{i}>\min F$.

Assume that the statement holds for $m$, and  let $F\in
\mc{S}_{m+1}$ and $(F_{i})_{i=1}^{d}$ be a partition of $F$ into
successive subsets.

We have also that  $F=\bigcup_{i=1}^{k}A_{i}$ where
$A_{i}\in\mc{S}_{m}$ and $k\leq A_{1}<A_{2}<\dots< A_{k}$. Let
$C_{i}=\{j\in \{2,3\dots, d\}: F_{j}\subset A_{i}\}$  for
$i=1,\dots, k$. We have that  $F_{j}\in\mc{S}_{m}$ for every $j$
with $F_{j}\subset A_{i}$ for some $i$. From inductive hypothesis
we get that
$$
U_{i}=\bigcup_{j\in C_{i}}G_{j}\in \mc{A}_{2}[\mc{S}_m]\,\,\,
\text{for every}\,\,  i\leq k.
$$
Let $U_{i}=U_{i,1}\cup U_{i,2}$ with $U_{i,1}<U_{i,2}$ in $\mc{S}_m$. Let us observe that
we have $3k\leq U_{i,1}$ for every $i=1,\dots,k$.

Consider also $B=\{j\leq d: F_{j}$ is not contained in any
$\mc{A}_{i}\}$.  It follows that $\# B\leq k$ and therefore
$\mc{L}=\{U_{i,1},U_{i,2}:\ i=1,\dots,k\}\cup\{G_{j}: j\in B\}$
contains at most $3k$ elements of $\mc{S}_m$ with $3k\leq\minsupp
G_{2}\leq\minsupp \mc{L}$ and therefore
$\cup\mc{L}\in\mc{S}_{m+1}$.

Since  $\mc{L}$ contains all $G_{j}$ except $G_{1}$ we get that
$$
\bigcup\{G_{j}:j\leq d\}\in\mc{A}_{2}[\mc{S}_{m+1}].
$$
The proof of the second part follows similarly, in a shorter way,
since we don't need the partition of sets $U_i$ into elements of
$\mc{S}_m$ and treating separately the set $G_1$.
\ep

\subsection{The quasiminimality}

We assume that in the whole section we work with a regular mixed Tsirelson space
$X=T[(\mc{S}_n,\theta_n)_{n\in\N}]$.

The arguments to prove the quasiminimality of the space $X$ will follow the arguments we
use for the $p-$spaces, showing the existence of special averages and forming from them
equivalent sequences.

We will need averages of a certain type, playing the role of
$\ell_p-$averages in $p-$spaces. The major tool providing these
averages is Lemma 4.11 \cite{ao} stated in general case, however
with a proof in \cite{ao} requiring an additional assumption -
that the sequence $(c_n)_n$ is decreasing. Thus we need this
assumption in the proof of Theorem \ref{t-quasi}, and we use it
only when applying Lemma 4.11 \cite{ao}. Therefore any mixed
Tsirelson space $T[(\mc{S}_n,\theta_n)_n]$ in which Lemma 4.11
holds, is quasiminimal.

Let us mention that as observed by the first named author and D. Leung (private
communication) the proof of the main result of \cite{ao} on distortion of spaces
$T[(\mc{S}_n,\theta_n)_n]$ with $c_n\to 0$ can be adapted to work in all cases just
following the proof of the main result from the preprint version of \cite{ao} in
arxiv.org.

In order to have an analogue to Lemma \ref{p-xi} we shall use the auxiliary space
$X_k=T[(\mc{S}_{n}[\mc{A}_k],\theta_n)_{n}]$, $k\in\N$. It follows readily from the
definition that $\nrm{\cdot}_{X}\leq\nrm{\cdot}_{X_k}$.  In the next lemma we prove that
the spaces are isomorphic.

\bl\label{l3.3} The spaces $X$ and  $X_k$ are $(k+1)$-isomorphic.

In particular for every functional $f$ in the norming set of $X_k$ there exists $(k+1)$
functionals $\phi_{1}<\dots<\phi_{k+1}$ in the norming set  of  $X$ such that
$f=\phi_{1}+\dots+\phi_{k+1}$ and the restriction of the tree-analysis of $f$ to the
support of $\phi_i$ gives a tree-analysis of $\phi_i$, $i=1,\dots, k+1$.

\el

\bp We prove it by induction on the height of the tree-analysis. If $f=\theta_n\sum_{i\in
A}e_i$, with $A\in \mc{S}_n[\mc{A}_k]$, then by Lemma \ref{sch1} $A\in
\mc{A}_{k+1}[\mc{S}_n]$ hence $A=A_1\cup \dots\cup A_{k+1}$ for some $A_t$'s from $\mc{S}_n$,
hence $\phi_t=\theta_n\sum_{i\in A_t}e_i$, $t=1,\dots,{k+1}$, give the desired partition.

Now take $f=\theta_n\sum_{i\in A}f_i$ for some $A\in \mc{S}_n[\mc{A}_k]$ and some norming
functionals $(f_i)$ of $X_k$. By inductive assumption
$$
f=\theta_n\sum_{i\in A}(\phi_{i,1}+\dots+\phi_{i,k+1})
$$
for some norming functionals $(\phi_{i,j})_{i\in A,j=1,\dots,k+1}$ of $X$.

Since $(f_{i})_{i\in A}$ is $\mc{S}_n[\mc{A}_k]$-admissible and $k(k+1)<2^{k+1}$, Lemma
\ref{sch1} yields $(\phi_{i,j})_{i\in A,j=1,\dots,k+1}$ is
 $\mc{A}_{k+1}[\mc{S}_n]$-admissible. Therefore there are some $G_1,\dots,G_{k+1}$ such that
$(\phi_{i,j})_{(i,j)\in G_t}$ is $\mc{S}_n$-admissible and setting
$\phi_t=\theta_n\sum_{(i,j) \in G_{t}}\phi_{i,j}$ for $t=1,\dots,k+1$ we get $(k+1)$
functionals $\phi_1<\dots<\phi_{k+1}$ in the norming set of $X$ with
$\phi_1+\dots+\phi_{k+1}=f$.

It follows that  $\abs{f(u)}\leq
(k+1)\max\{\abs{\phi_{1}(u)},\dots,\abs{\phi_{k+1}(u)}\}$ for any $u\in X$. The
restriction of the tree-analysis of $f$ to the support of each $\phi_t$ gives the
tree-analysis of $\phi_t$ in $X$ with the desired property. \ep

We prove now an analogue to Lemma \ref{p-xi}.

\bl \label{t-xi} Take a functional $\phi$ from the norming set of
$X$ and any finite block sequence $(v_n)$. Set $v=\sum_n v_n$.
Then there is a functional $\phi'$ with the tree-analysis
$(\phi'_t)_{t\in T'}$ comparable with $(v_n)$ and $4\phi'(v)\geq
\phi (v)$.
\el

\bp For the proof we consider the auxiliary space $X_3$. Take $\phi$ in the norming set
of $X$ with a tree-analysis $(\phi_t)_{t\in T}$. For any $t\in T$ set $I_{t}=\{n: v_{n}$
is covered by $\phi_{t}\}$ and set for every $n\in I_{t}$, $J_{t,n}=\{r\in \suc t:
\ran\phi_{r}\cap\ran v_{n}\ne\emptyset\}$.

Observe that every $r\in \suc t$ belong to at most two sets $J_{t,n}$, say
$J_{t,n_{1}^r}, J_{t,n_{2}^r}$. So we can split each $\supp\phi_{r}$ into three norming
functionals of $X$ with a tree-analysis comparable to $(v_n)$: $(\phi_r)|_{\supp
v_{n_1^r}}$, $(\phi_r)|_{\supp v_{n_2^r}}$ and the remaining part.

Applying the above argument for all $t\in T$ with $I_{t}\ne \emptyset$ and $r\in\suc t$, starting from the root and moving ``downward'' the tree,
we produce a tree-analysis comparable with $(v_n)$ of some functional $f$ in the norming
set of $X_3$.

By Lemma  \ref{l3.3} the functional $f$ of $X_3$ can be
represented as the sum of four norming functionals
$\phi_1<\phi_2<\phi_3<\phi_4$ of $X$ such that restriction of the
tree-analysis of $f$ to the support of each $\phi_l$ gives the
tree-analysis of $\phi_l$. Therefore each $\phi_l$ has a
tree-analysis comparable with $(v_n)$. To end the proof pick
$\phi'$ from $\{\phi_1,\dots,\phi_4\}$ satisfying $\phi'(v)\geq
\phi(v)/4$. \ep

We recall now the notion of averages given in \cite{ao},
Definition 4.5, with simplified notation. Averages of this type,
called special convex combinations, were also used in studies of
mixed Tsirelson spaces in particular in \cite{ad,adm,lt,klmt}.

\bd \cite{ao} Fix $M\in\N$, $\e>0$ and take a block sequence
$(x_n)\subset X$. We say that $x\in X$ is an $(M,\e)$
\textit{average} of $(x_n)$ given by an \textit{averaging tree}
$T=(x^{i,j})_{j=0,i=1}^{M,N_j}$ on $X$ if

\bnum

\item $1=N_M<N_{M-1}<\dots<N_0$,

\item $x=x^{1,M}$ is the root of $T$, the sequence $(x^{i,0})_{i=1}^{N_0}$ of maximal nodes
of $T$ form a subsequence of $(x_n)$,

\item $\suc x^{i,j}=(x^{s,j-1})_{s\in I_{i,j}}$ for some non-empty interval
$I_{i,j}\subset\{1,\dots,N_j\}$ with $(x^{s,j-1})_{s\in I_{i,j}}$ a block sequence
$\mc{S}_1$-admissible with respect to $(x_n)$, for any $j=0,\dots,M$, $i=1,\dots,N_j$,
\item $x^{i,j}=\frac{1}{k_{i,j}}\sum_{s\in I_{i,j}}x^{s,j-1}$, where $k_{i,j}=\#
I_{i,j}$,  for any $j=0,\dots,M$, $i=1,\dots,N_j$,

\item $k_{1,j}>6\cdot 2^{2+j}\theta^{-1}\e^{-1}$, $k_{i,j}>6\cdot
2^{1+i+j}\theta^{-1}\e^{-1}\max\supp x^{i-1,j}$ for any $i=2,\dots, N_j$, $j=0,\dots,M$.
\enum \ed

\br Recall that by Proposition 4.7, \cite{ao}, for any normalized
block sequence $(x_n)$, $M\in\N$, $\e>0$ there is a
$(M,\e)-$average of $(x_n)$.
\er

\bl\label{theta} For any block subspace $Y$ of $X$, any $\de>0$, $\e>0$ and any $M\in\N$
there is a normalized block sequence $(x_i)\subset Y$ and a $(M,\e)-$average $x$ of
$(x_i)$ given by an averaging tree $T=(x^{i,j})_{j=0,i=1}^{M,N_j}$ such that for any
$j=1,\dots,M$ we have
$$
\nrm{x^{i,j}}\geq (1-\de)^j\theta^j\ \ \ i=1,\dots,N_j
$$
\el

\bp  We use the standard argument as in Lemma 2.8 \cite{ad},   Lemma 4.12 \cite{ao}.
Assume the contrary and pick any normalized block sequence $(x_n)$ in $Y$ and let
$(y_{n,1})$ be a block sequence of $(M,\e)-$averages of $(x_n)$. Let any $y_{n,1}$ be
given by a tree $(x^{i,j}_{n,1})_{i,j}$. For any $n$ there is some $1\leq j_{n,1}\leq M$
and $i_{n,1}$ so that $\nrm{x^{i_{n,1},j_{n,1}}_{n,1}}\leq
(1-\de)^{j_{n,1}}\theta^{j_{n,1}}$. Hence there is an infinite $N_1\subset\N$ and some
$1\leq J_1\leq M$ so that $j_{n,1}=J_1$ for any $n\in N_1$.

Let $z_n=x^{i_{n,1},J_1}_{n,1}/\nrm{x^{i_{n,1},J_1}_{n,1}}$ for any $n\in N_1$. Let
$(y_{n,2})$ be a block sequence of $(M,\e)-$averages of $(z_n)_{n\in N_1}$, each
$y_{n,2}$ given by an averaging tree $(x^{i,j}_{n,2})_{i,j}$. Find some infinite
$N_2\subset N_1$, some sequence of integers $(i_{n,2})$ and some $1\leq J_2\leq M$ so
that $\nrm{x^{i_{n,2},J_2}_{n,2}}\leq (1-\de)^{J_2}\theta^{J_2}$ for any $n\in N_2$.
Notice that by definition of averages each $x^{i_{n,2},J_2}$ is a convex combination of
some $\mc{S}_{J_2}$-admissible $z_n$'s, and recall that any $x^{i,j}_{n,1}$ is a convex
combination of some $\mc{S}_{J_1}$-admissible $x_n$'s, so
$x_{n,2}^{i_{n,2},J_2}=\sum_{k\in F_{n,2}}a_kx_k$, for $(x_k)_{k\in F_{n,2}}$
$\mc{S}_{J_1+J_2}$ -admissible and $\sum_{k\in F_{n,2}}a_k\geq
(1-\de)^{-J_1}\theta^{-J_1}$. Normalize $(x^{i_{n,2},J_2})$ and continue in this manner.

After $m$ steps we end up with a vector $x=x^{i_{1,m},J_m}_{1,m}$
(first of some sequence) so that
$\nrm{x}=\nrm{x^{i_{1,m},J_m}_{1,m}}\leq
(1-\de)^{J_m}\theta^{J_m}$. On the other hand, $x=\sum_{k\in
F_{1,m}}a_kx_k$, for some $(x_k)_{k\in F_{1,m}}$
$\mc{S}_{J_1+\dots+J_m}$-admissible and
$$
\sum_{k\in F_{1,m}}a_k\geq (1-\de)^{-(J_1+\dots+J_{m-1})}\theta^{-(J_1+\dots+J_{m-1})}
$$
Hence we get
$$
(1-\de)^{J_m}\theta^{J_m}\geq \theta_{J_1+\dots+J_m} (1-\de)^{-(J_1+\dots+J_{m-1})}\theta^{-(J_1+\dots+J_{m-1})}
$$
which yields $\theta_{J_1+\dots+J_m}\leq
(1-\de)^{J_1+\dots+J_m}\theta^{J_1+\dots+J_m}$, which since
$J_i\geq 1$, for large $m$ contradicts the definition of $\theta$.
\ep

Now we are ready to prove the existence of special $(M,\e)-$averages in all block
subspaces of a mixed Tsirelson space. Recall that this fact is well-known and used in
case of $\theta=1$, e.g. \cite{ad,adm,m,lt,klmt}.

\bc\label{t-av} Let $X=T[(\mc{S}_n,\theta_n)_n]$ be a regular
mixed Tsirelson space with the sequence
$(\frac{\theta_n}{\theta^n})_n$ decreasing. Then for any block
subspace $Y$ of $X$ and any $M\in\N$ there is a vector $y\in Y$
with $\nrm{y}=1$ such that for any $0\leq j\leq M$
$$
4\theta_1^{-1}\theta^{-j-1}\geq\sup\left\{\sum_i\nrm{E_iy}:\ (E_i)
\ - \ \mc{S}_j\ \mathrm{admissible}\right\}\geq
\theta_1\theta^{1-j}/4
$$
\ec

\bp Pick the $(M,\e)-$average $x$ of $(x_n)$ provided  by Lemma \ref{theta}, for $\de$
with $(1-\de)^M\geq 1/2$ and $\e\leq \theta^{M-1}\theta_M$. It is enough to show that
$$
2\theta^{M-1}\geq \nrm{x}\geq \theta_1\theta^{M-1}/2
$$
$$
2\theta^{M-j-1}\geq\sup\left\{\sum_i\nrm{E_ix}:\ (E_i) \ - \
\mc{S}_j\ \mathrm{admissible}\right\}\geq \theta^{M-j}/2, \ \
j=1,\dots,M
$$
The lower estimate for $j=1,\dots,M$ is provided by Lemma
\ref{theta}. Let $x$ be given by a tree $(x^{i,j})_{i,j}$. Notice
that for any $j=1,\dots,M-1$ we have
$x=\sum_{i=1}^{N_{M-j}}a_i^jx^{i,M-j}$ for some $(a_i^j)\subset
[0,1]$ with $\sum_{i=1}^{N_{M-j}}a_i^j=1$ and
$(x^{i,M-j})_{i=1}^{N_{M-j}}$ is $\mc{S}_j$-admissible, which ends
the proof. For $j=0$ it is enough to use estimation for $j=1$ and
definition of the norm.

The upper estimate for any $j=1,\dots,M$ follows from the proof of
Lemma 4.11 (7), \cite{ao}, in case $(\frac{\theta_n}{\theta^n})_n$
decreasing, for $N=1$. The upper estimate for $j=0$, i.e. the
estimate of $\nrm{x}$ follows again from the proof Lemma 4.11 (6),
\cite{ao}, in case $(\frac{\theta_n}{\theta^n})_n$ decreasing, for
$N=1$, $J=M$, $i=1$. \ep

Thanks to more complex structure of families $(\mc{S}_n)_{n\in\N}$, which can "absorb"
families $\mc{A}_n$ as it is reflected in Lemma \ref{sch1}, the construction of
equivalent sequences in the spaces $T[(\mc{S}_n,\theta_n)_{n\in\N}]$ is simpler than the
one in case of $p-$spaces. Namely in $p-$spaces case we produce equivalent sequences from
averages of $\ell_p-$averages. Now the equivalent sequence consist simply of ''special
averages'' provided by Corollary \ref{t-av}.

\bt\label{t-quasi} Any regular mixed Tsirelson space
$X=T[(\mc{S}_n,\theta_n)_{n\in\N}]$, with the sequence
$(\frac{\theta_n}{\theta^n})_n$ decreasing, is quasiminimal. \et

By Corollary \ref{t-av} in order to prove the Theorem it is enough
to show the following

\bl\label{t-eq} Let $(y_n), (z_n)\subset X$ be normalized block sequences satisfying for
some $(m_n)\subset\N$ the following:

\bnum

\item $2\maxsupp y_{n}<z_{n}$ and $3\maxsupp z_{n}<y_{n+1}$ for any $n\in\N$,

\item $y_{n}$ and $z_{n}$  satisfies the conclusion of Corollary \ref{t-av} with $M=m_n$ for every $n\in\N$,

\item $\theta_{m_{n+1}}\maxsupp z_n<2^{-n}$, $\theta_{m_{n+1}}\maxsupp y_{n}<2^{-n}$ for
any $n\in\N$.

\enum

Then $(y_n)_n$ and $(z_n)_n$ are equivalent.

\el

\bp[Proof of Lemma \ref{t-eq}]

We show first that $(z_{n})_{n}$ dominates $(y_{n})_{n}$ with constant $C=(4\cdot 2)\cdot
(2\cdot \frac{4^2}{\theta_{1}^2\theta^3})$. Take any scalars $(a_n)_{n\in D}$. We shall
use again the auxiliary space $X_{2}=T[(\mc{S}_{n}[A_{2}], \theta_{n})_{n}$. Namely for
every functional $\phi$ in the norming set of $X$ with a tree-analysis $(\phi_{t})_{t\in
T}$ comparable with $(y_n)$ we shall produce two functionals $f_{1},f_{2}$ in the norming
set of $X_{2}$ such that
\begin{equation}\label{ee1}
\phi(\sum_{n\in D}a_{n}y_{n})\leq (2\cdot
\frac{4^2}{\theta_{1}^2\theta^3})(f_{1}+f_{2})(\sum_{n\in D}a_{n}z_{n})
\end{equation}
By Lemma \ref{l3.3} and Lemma \ref{t-xi} this will prove that $\nrm{\sum_na_ny_n}\leq
C\nrm{\sum_na_nz_n}$.

For every $t\in T$ we set
$$
D_{t}=\{n\in\N: \supp \phi_{t}\cap\supp y_{n}=\supp \phi\cap\supp y_{n}\}
$$
and  as usual  $I_{t}=\{n\in D_{t}: \phi_{t}\,\,\text{covers}\,\,  y_{n}\}$.

We shall produce inductively (on level of $t\in T$) for any $t\in T$ two functionals
$f_{t,1}<f_{t,2}$ in the norming set of $X_2$ with supports contained in $\bigcup\{\supp
z_n:\ n\in D_t\}$ such that \eqref{ee1} holds for $\phi_{t}, D_{t}$ and $f_{t,1},
f_{t,2}$.

Assume $t$ is maximal, then $\phi_t=e^*_i$ for some $n_0$. If $D_t=\emptyset$ let
$f_{t,1}=f_{t,2}=0$, otherwise $D_t=I_{t}=\{n\}$ and let $f_{t,1}$ be the norming
functional of $z_n$ in $X$, $f_{t,2}=0$.

Assume that we have chosen suitable pairs of functionals for any functionals on the
levels $l+1,l+2,\dots\text{height}(T)$ and take $\phi_{t}$ on level $l$. Let
$\phi_{t}=\theta_{k}\sum_{s\in \suc t}\phi_{s}$. As usual for $n\in I_{t}$ we set
$J_{n}=\{s\in \suc t: \supp \phi_{s}\subset \ran y_{n}\}$.

\

CASE 1.  $k\leq m_{n}$ for every $n\in I_{t}$.

Then if $(\phi_{s})_{s\in J_{n}}$ is $\mc{S}_{k_{n}}$- and not
$\mc{S}_{k_{n}-1}$-admissible by condition 2. in choice of $(y_n)$ we get
$$
\sum_{s\in J_{n}}\phi_{s}(y_{n})\leq   4\theta_{1}^{-1}\theta^{-k_n-1}
$$
Hence by condition 2. in choice of $(z_n)$ we can take an $\mc{S}_{k_n-1}$-admissible set
of functionals $(\psi_{u})_{u\in V_{n}}$ from the norming set of $X$ with $\supp
\psi_u\subset\supp z_n$, $u\in V_n$, such that
$$
\sum_{s\in J_{n}}\phi_{s}(y_{n})\leq \frac{4^2}{\theta_{1}^{2}\theta^3}\sum_{u\in
V_{n}}\psi_{u}(z_{n})
$$
>From the inductive hypothesis for every $s\in \suc t$ such that
$D_{s}\ne\emptyset$ we have take two suitable functionals
$f_{s,1}<f_{s,2}$ with $\supp f_{s,1},\supp f_{s,2}\subset\{\supp
z_{n}:n\in D_{s}\}.$ From Lemma \ref{l3}(b) we get that the set
$$\{\supp\psi_{u}: u\in V_{n}, n\in I_{t}\}\cup\{\supp f_{s,1}: D_s\ne\emptyset\} $$
is $\mc{S}_{k}$-admissible and therefore the sequence
$$
\bigcup_{n\in I_{t}}\{\psi_{u}: u\in V_{n}\}\cup \bigcup_{s\in\suc
t:D_{s}\ne\emptyset}\{f_{s,1},f_{s,2}\}
$$
is $\mc{S}_k[\mc{A}_2]$-admissible. It follows that
$$
f_{t,1}=\theta_{k}\left(\sum_{n}\sum_{u\in
V_{n}}\psi_{u}+\sum_{s:D_{\phi_{s}\ne\emptyset}}(f_{s,1}+f_{s,2})\right)
$$
is in the norming set of the set of the space $X_{2}$ and
$$
\phi_{t}(\sum_{n\in D_{t}}y_{n})\leq \frac{4^2}{\theta_{1}^2\theta^3}f_{t,1}(\sum_{n\in
D_{t}}a_nz_n).
$$
In this case we set $f_{t,2}=0$.

\

CASE 2. There exists  $n\in I_{t}$ such that  $m_{n}<k$.

Then let  $n_0$ be the unique $n\in I_{t}$ such that  $m_{n}\leq
k<m_{n+1}$. For every $n\in D_{t}$ with $ n<n_{0}$ we have
$\theta_{k}\maxsupp (y_{n})\leq 2^{-n}$ so with error $\max\vert
a_n\vert\sum_{n\in D_{t}, n<n_0}2^{-n}$ we can  erase this part.

As usual we compare $\phi_t(a_{n_0}y_{n_0})$ and  $\phi_t(\sum_{I_{t}\ni n>n_0}a_ny_n)$.
If the second  term  dominates the first one  we proceed as in Case 1, but multiplying by
2 the estimation of $\phi_t(\sum_{n\in I_{t}}a_ny_n)$.

If the first term dominates the second one, take the functional
$$
f_{t,2}=\theta_{k}\sum_{s\in\suc t:D_s\ne\emptyset}(f_{s,1}+f_{s,2}).
$$
Since the set $\{\phi_{s}: D_{s}\ne\emptyset\}$ is $\mc{S}_{k}-$admissible, it is clear
that $f_{t,2}$ is in the norming set of $X_{2}$. Now take a functional $f_{t,1}$ that
norms $z_{n_0}$ in $X$. Of course we may assume that $\ran f_{t,1}\subset\ran z_{n_0}$.
We also have that $f_{t,1}<f_{t,2}$ since we have delete all $z_{n}$ for $n\in D_{t}$ and
$n<n_0$. Hence we have produce two functionals $f_{t,1}<f_{t,2}$ in the norming set of
$X_{2}$ satisfying the desired conditions.

\

We prove now that  $(y_{n})_{n}$ dominates $(z_{n})_{n}$ in an analogous way.

Fix scalars $(a_n)_{n\in D}$. Take any norming functional $\phi$ and assume by Lemma
\ref{t-xi} that $\phi$ has a tree-analysis $(\phi_{t})_{t\in T}$ comparable with $(z_n)$.

For every $t\in T$ with $D_{t}=\{n: \supp \phi_{t}\cap\supp z_{n}=\supp \phi\cap\supp
z_{n}\}\neq \emptyset$ we shall produce two functionals $f_{t,1}<f_{t_2}$ in the norming
set of the space $X_{2}$ such that $\supp f_{t,1},\supp f_{t,2}\subset\bigcup\{\supp y_n:
\ n\in D_t\}$ and
$$
\phi_{t}(\sum_{n\in D_{t}}a_{n}z_{n})\leq C(f_{t,1}+f_{t,2})\left(\sum_{n\in
D_{t}}a_{n}y_{n}\right)
$$
If $t$ is maximal we proceed as in the previous part. Assume we have done it for all
$t\in T$ in levels $l+1,\dots,height(T)$. Let $\phi_{t}=\theta_{k}\sum_{s\in \suc
t}\phi_{s}$, for some $\{\phi_{s}:\ s\in \suc t\}-\mc{S}_{k}$-admissible.

We set
$$
I_{t}=\{n\in D_{t}: \phi_{t}\,\,\text{covers}\,\,z_{n}\}\,\,\text{and}\,\,\,B=\{s\in \suc t:\ D_{s}\ne\emptyset\}.
$$
For every $n\in I_{t}$ we set  $J_{n}=\{s\in \suc t: \supp \phi_{s}\subset\ran y_n\}$.

\

CASE 1. $k\leq  m_{n}$ for every $n\in I_{t}$.

We have that $\suc t=\bigcup_{n\in I_{t}}\{s: s\in J_{n}\}\cup \bigcup_{s\in B}\{s\}.$ By
condition 2. of choice of $(y_n)$ and $(z_n)$ for every $n\in I_{t}$ we choose
$\{\psi_{u}: u\in V_{n}\}$ and $\mc{S}_{k_{n}-1}$-admissible set, where  $\mc{S}_{k_n}$
is the admissibility of the set$\{\phi_{s}: s\in J_{n}\}$, such that
$$
\sum_{s\in J_{n}}\phi_{s}(z_n)\leq \frac{4^2}{\theta_{1}^2\theta^3}\sum_{u\in
V_{n}}\psi_{u}(y_{n}).
$$
By Lemma \ref{l3} we get that the set
$$\bigcup_{n\in I_{t}}\{\psi_{u}: u\in V_{n}\}\cup\{f_{s,1}: s\in B\}
$$
is $\mc{A}_{2}[\mc{S}_{k}]$-admissible. Therefore we have
partition the set
\begin{align*}
\bigcup_{n\in I_{t}}\{& \psi_{u}: u\in V_{n}\}\cup\{f_{s,1}: s\in B\}\\
&=(\{\psi_u:\ u\in
V_1\}\cup\{f_{s,1}:\ s\in B_1\})\cup (\{\psi_u:\ u\in V_2\}\cup\{f_{s,1}:\ s\in B_2\})
\end{align*}
into two successive $\mc{S}_k$-admissible sequences. Therefore the
functionals
$$
f_{t,i}=\theta_k\left(\sum_{u\in V_i}\psi_u+\sum_{s\in B_i}(f_{s,1}+f_{s,2})\right), \ \
\ i =1,2
$$
are from the norming set of $X_2$, satisfying the desired estimation. Indeed
\begin{align*}
\phi_{t}(\sum_{n\in D_{t}}a_{n}z_{n})&=\theta_{k}\sum_{n\in I_{t}}\sum_{s\in
J_{n}}\phi_{s}(a_nz_n)+\theta_{k}\sum_{s\in B}\phi_{s}(\sum_{n\in D_{s}}a_{n}z_{n})
\\
&\leq \theta_{k}\sum_{n\in I_{t}}\sum_{u\in V_{n}}C\psi_{u}(a_ny_n) +\theta_{k}\sum_{s\in
B}C(f_{s,1}+f_{s,2})(\sum_{n\in D_{s}}a_{n}y_{n})
\\
&=C(f_{t,1}+f_{t,2})\left(\sum_{n\in D_{t}}a_{n}y_{n}\right).
\end{align*}

\

CASE 2. There exist  $n\in I_{t}$ such that  $m_{n}\leq k$.

Take $n_0=\max\{n: m_{n}\leq k\}$.  By condition 3. in the choice of the sequences
$(y_n),(z_n)$ we can erase the $y_{n}$, $n<n_0, n\in D_{t}$ with error $\sum_{n\in D_{t},
n<n_0}\vert a_{n}\vert 2^{-n}$.

Now we compare $\phi_{t}(a_{n_0}z_{n_0})$ and $\theta_{k}\sum_{I_{t}\ni n>n_0}\sum_{s\in
J_{n}}\phi_{s}(a_{n}z_{n})$. Assume  that the first term dominates the second. It follows
$$
\phi_{t}(a_{n_0}z_{n_0})+\theta_{k}\sum_{I_{t}\ni n>n_0}\sum_{s\in
J_{n}}\phi_{s}(a_{n}z_{n})\leq 2\phi_{t}(a_{n_0}z_{n_0}).
$$
In this case we take a functional $f_{t,1}$ that norms $y_{n_0}$
and we set
$$
f_{t,2}=\theta_{k}\sum_{s\in B}(f_{s,1}+f_{s,2}).
$$
It is clear that  $f_{t,1}<f_{t,2}$ and by Lemma \ref{l3} we get that $\{f_{s,1}: s\in
B\}$ is  $\mc{S}_{k}$-admissible hence $f_{t,2}$ is in the norming set of $X_{2}$.
Therefore
\begin{align*}
\phi_{t}(\sum_{n\in D_{t}}a_{n}z_{n})&=\phi_{t}(a_{n_0}z_{n_0})+\theta_{k}\sum_{n_0<n\in
I_{t}}\sum_{s\in J_{n}}\phi_{s}(a_nz_n)+\theta_{k}\sum_{s\in B}\phi_{s}(\sum_{n\in
D_{s}}a_{n}z_{n})
\\
&\leq 2\phi_{t}(a_{n_0}z_{n_0})+\theta_{k}\sum_{s\in B}C(f_{s,1}+f_{s,2})(\sum_{n\in
D_{s}}a_{n}y_{n})
\\
&\leq Cf_{t,1}(a_{n_0}y_{n_0})+\theta_{k}\sum_{s\in B}C(f_{s,1}+f_{s,2})(\sum_{n\in
D_{s}}a_{n}y_{n})
\\
&=C(f_{t,1}+f_{t,2})\left(\sum_{n\in D_{t}}a_{n}y_{n}\right).
\end{align*}
If the second term dominates the first one we repeat the argument
by multiplying the estimation on $\phi_{t}(\sum_{n\in
I_{t}}a_nz_n)$ by $2$. \ep

\br The reasoning presented in the proof of Theorem \ref{t-quasi} can be adapted in case
of $\theta=1$ to obtain a version of Proposition 5.7 \cite{m}, giving equivalence of
$(y_n)$ of long "special" averages to some subsequence of the basis - it is enough to use
the fact that in this case unit basis vectors $(e_n)$ satisfies the assertion of
Corollary \ref{t-av}.

\er

\br Notice that as in the case of $p-$spaces the above proof provides in every block
subspace of $X$ some subsequentially minimal block subspace. Indeed any block sequence of
vectors of the form $(y_1+\dots+y_N)/\nrm{y_1+\dots+y_N}$, where $y_1,\dots,y_N$ are
chosen as in the beginning of the proof of Theorem \ref{t-quasi} for arbitrary big $N$
and $m_1,\dots,m_N$, spans a subsequentially minimal subspace.

We recall here that by Theorem 35 \cite{klmt} if $0<\inf
c_n\leq\sup c_n<1$ then the space $X=T[(\mc{S}_n,\theta_n)_n]$
fails to be a subsequentially minimal space in a strong way: $X$
is saturated with block subspaces in which no block sequence is
equivalent to a subsequence of the basis. It was also shown that
the assumption $\sup_n c_n<1$ is not necessary to have the
property described above.

\er

\section{A note for the quasiminimality of the dual spaces}

In this section we give a general result which allows to transfer
the minimality and quasiminimality properties from a space to its
dual space.

\bt\label{t1} Let $X$ be a Banach space  with shrinking basis
$(e_{n})_{n\in\N}$ and let $(x_{n}^{*})_{n\in\N}$ be a block
sequence  of $(e_{n}^{*})_{n}$ such that there  exist a block
sequence $(x_{n})_{n}$ of the basis such that
\begin{enumerate}
\item $\frac{1}{1+\e}\leq\nrm{x_n}\leq 1+\e$ and
$x_{n}^{*}(x_n)=1$.

\item  $m\nrm{\sum_{n}a_ne_{k_n}}\leq \nrm{\sum_na_nx_n}\leq
M\nrm{\sum_na_ne_{k_n}}$ for some  $k_n\to\infty$.

\item The map  $P:X\to [x_n]$ defined by  $P(x)=\sum_{n}x_{n}^{*}(x)x_{n}$ is a bounded
projection.
\end{enumerate}
Then   $(x_n^*)_{n}$ is equivalent to the subsequence
$(e_{k_n}^{*})$. \et

\bp From the inequality $\nrm{\sum_na_nx_n}\leq
M\nrm{\sum_na_ne_{k_n}}$ it follows that
$$
\nrm{\sum_nc_ne_{k_n}^{*}}\leq M\nrm{\sum_nc_nx_{n}^{*}}
$$
For the converse inequality let $x\in \mc{S}_{X}$ be such that
$\sum_{n}c_nx_{n}^{*}(x)\geq
(1+\e)^{-1}\nrm{\sum_{n}c_{n}x_{n}^{*}}$. We have that
 $$
 \nrm{\sum_{n}x_{n}^{*}(x)e_{k_n}}\leq
 \nrm{\sum_{n}x_{n}^{*}(x)x_{n}}
 =\nrm{P(x)}\leq \nrm{P}
 $$
It follows that
\begin{align*}
\nrm{\sum_{n}c_ne^{*}_{k_n}}&\geq
\nrm{P}^{-1}\sum_{n}c_{n}e^{*}_{k_n}(\sum_{m}x_{m}^{*}(x)e_{k_m})
\\
&= \nrm{P}^{-1}\sum_{n}c_{n}x_{n}^{*}(x) \geq ((1+\e)\nrm{P}
)^{-1}\nrm{\sum_{n}c_{n}x_{n}^{*}}
\end{align*}
\ep

\bc The dual space of Schlumprecht space
$S=T[(\mc{A}_n,1/\log_2(n+1))_n]$ is minimal (\cite{CKKM}). Dual
spaces of mixed Tsirelson spaces $T[(\mc{S}_n, \theta_n)_n]$ with
$\sup\theta_n^{1/n}=1$ are quasiminimal.
\ec

\bp In order to prove the Corollary we recall the following

\bd Let $X$ be a Banach space with a basis. A normalized vector
$x$ is called a $(1+\e)-c_{0}^{N}-$average, if
$x=\sum_{k=1}^{N}x_{k}$ for some finite block sequence
$(x_{k})_{k=1}^{N}$ with $(1+\e)^{-1}\leq \nrm{x_k}\leq 1$ for all
$k$. \ed
\br\label{c0aver} Assume $X$ is a Banach space with a
shrinking unconditional basis. With any $(1+\e)-c_0^N-$average
$x^*=\sum_{k=1}^Nx_k^*\in X^*$ we can associate
$(1+\e)-\ell_{1}^{N}-$average in $X$ in the following way: for
every $k$ pick a normalized $x_k\in X$ be such that $\supp
x_k=\supp x_k^*$  and $x_{k}^{*}(x_k)=\nrm{x_k^{*}} \geq
(1+\e)^{-1}$. Then the vector
$x=\sum_{k=1}^{N}x_{k}/\nrm{\sum_{k=1}^{N}x_{k}}$ is
$(1+\e)-\ell_1^N-$average, since for any non-negative scalars
$(a_k)_{k=1}^N$ we have
$$
\nrm{\sum_ka_kx_k}\geq
x^{*}(\sum_ka_kx_k)=\sum_{k=1}^{N}a_kx_{k}^{*}(x_k)\geq
(1+\e)^{-1}\sum_ka_k.
$$
\er

\bd\cite{at} Let $X$ be a Banach space with the basis
$(e_k)_{k\in\N}$. Let $\e> 0$ and $j\in\N$, $j > 1$. A convex
combination $\sum_{k\in F} a_ke_k$ of the basis $(e_k)_{k\in\N}$
is called an $(\e; j)$ \textit{basic special convex combination}
if $F\in \mc{S}_j$ and $\sum_{k\in G} a_k < \e$ for every $G\in
\mc{S}_{j-1}$.

Let $(x_k)_{k\in\N}$ be a block sequence. A convex combination $
\sum_{k\in F} a_kx_k$ of the sequence $(x_k)_{k\in\N}$ is called
an $(\e; j)$ \textit{special convex combination} ($(\e; j)$
s.c.c.) of $(x_k)_{k\in\N}$ if $\sum_{k\in F} a_ke_{t_k}$ (where
$t_k = \min\supp x_k$ for each $k$) is an $(\e; j)$ basic special
convex combination.

\ed

In order to apply Theorem \ref{t1} we need the following result.
\bl\label{duallemma} Let $\e\in (0,1)$ and $Y^*$ be a subspace of
$S^*$, the dual space to Schlumprecht space. Then for every
$N\in\N$ there exists a $(1+\e)-c_0^N-$average $z^*\in Y^*$.

Let $\e\in (0,1)$ and $Y^*$ be a subspace of the dual space of a
mixed Tsirelson space $X=T[(\mc{S}_n,\theta_n)_n]$ with
$\sup\theta_n^{1/n}=1$. Then for every $N\in\N$ there exists an
$(\e,N)$ s.c.c. $z=\sum_{i\in F}a_iz_i\in X$ and $z^*\in Y^*$ with
$z^*(z)\geq 1-\e$.
\el

The above Lemma for mixed Tsirelson spaces $T[(\mc{S}_n,
\theta_n)_{n}]$ has been proved in \cite{at}, Proposition 5.4. For
Schlumprecht space this lemma appeared in \cite{CKKM}, Prop. 5.

\

We are now ready to complete the proof of the corollary.

Let $(z_{k}^{*})$ be a block sequence of  $2-c_0^{N_k}-$averages in the block subspace
$Y^*$ of $S^*$ and  $(z_{k})_k$ be the sequence of the associated $\ell_1^{N_k}-$averages
(Remark \ref{c0aver}). We assume that the sequence $N_k$ is fast increasing. By
\cite{as}, see also \cite{m}, the sequence $(z_{k})_{k}$ is equivalent to the basis of
$S$ and moreover spans a complemented subspace of $S$. Theorem \ref{t1} implies the
result for the space $S$.

For the mixed Tsirelson spaces the analogous properties of a rapidly increasing sequence
of s.c.c. averages have been proved for example in \cite{m}, thus again application of
Theorem \ref{t1} finishes the proof of Corollary. \ep

\br The dual spaces of the spaces $T[\mathcal{S}_{\alpha},
\theta]$ are minimal. Indeed, we use the following result:
\begin{proposition}\label{lm}\cite{LM}
Let normalized block sequences $(x_n)$, $(y_n)\subset
T[\mathcal{S}_{\alpha},\theta]$ be such that $x_n<y_n <x_{n+1}$
($n\in \N$). Then $(x_n)$ and $(y_n)$ are
$24\theta^{-2}$-equivalent.
\end{proposition}
Standard dual argument yields that  the same property holds also
in the dual space i.e. if any normalized block sequences
$(x_n^*)_n, (y_n^*)_n\subset T[\mathcal{S}_{\al},\theta]$ with
$x_n^*<y_n^* <x^*_{n+1}$,$n\in \N$ are equivalent. Now since $c_0$
is block finitely represented in the dual space
$T[\mc{S_\alpha},\theta]^*$ we follow the proof of \cite{CO} for
Tsirelson space $T[\mc{S}_1,\theta]$ to obtain the minimality.
\er

\end{document}